%% file: main.tex
\documentclass[11pt]{article}
\usepackage{algorithm}
\usepackage{algpseudocode}
\usepackage{algorithmicx}
\usepackage{amsfonts}
\usepackage{amsmath}
\usepackage{amssymb}
\usepackage{amsthm}
\usepackage{array}
\usepackage{bm}
\usepackage{cancel}
\usepackage{caption}
\usepackage{color}
\usepackage{colortbl}
\usepackage{diagbox}
\usepackage{enumitem}
\usepackage{fancyhdr}
\usepackage{float}
\usepackage{fullpage}
\usepackage{graphicx}
\usepackage{hyperref}
\usepackage{lastpage}
\usepackage{makecell}
\usepackage{mathtools}
\usepackage{multicol}
\usepackage{parskip}
\usepackage{setspace}
\usepackage{soul}
\usepackage{tasks}
\usepackage[normalem]{ulem}
\usepackage{url}
\usepackage{xparse}
\usepackage{lipsum}
\usepackage[svgnames]{xcolor}
\usepackage{tikz-cd}
\usepackage{pgfplots}
\usepackage{titling}
\usepackage[backend=bibtex, style=numeric]{biblatex}

\usepgfplotslibrary{fillbetween}

\pdfmapfile{-mpfonts.map} 

\makeatletter
\def\th@plain{
  \thm@notefont{}
  \itshape 
}
\def\th@definition{
  \thm@notefont{}
  \normalfont 
}
\makeatother

\newtheorem{thm}{Theorem}[section]
\newtheorem*{thm*}{Theorem}

\newtheorem{lem}[thm]{Lemma}

\newtheorem{cor}[thm]{Corollary}

\theoremstyle{definition}
\newtheorem{defn}[thm]{Definition}
\newtheorem{note}[thm]{Note}
\newtheorem{obs}[thm]{Observation}
\newtheorem{remark}[thm]{Remark}
\newtheorem{ex}[thm]{Exemple}
\newtheorem*{remark*}{Remark}

\newenvironment{Proof}[1][Proof]
{\proof[#1]\leftskip=0.5cm}
{\endproof}

\newenvironment{Sketch}[1][Sketch]
{\proof[#1]\leftskip=0.5cm}
{\endproof}

\usetikzlibrary{decorations.markings}
\usetikzlibrary{shapes.geometric}

\tikzstyle{none}=[inner sep=0pt]
\tikzstyle{rn}=[circle,fill=Red,draw=Black,line width=0.8 pt]
\tikzstyle{gn}=[circle,fill=Lime,draw=Black,line width=0.8 pt]
\tikzstyle{yn}=[circle,fill=Yellow,draw=Black,line width=0.8 pt]
\tikzstyle{simple}=[-,draw=Black,line width=2.000]
\tikzstyle{arrow}=[-,draw=Black,postaction={decorate},decoration={markings,mark=at position .5 with {\arrow{>}}},line width=2.000]
\tikzstyle{tick}=[-,draw=Black,postaction={decorate},decoration={markings,mark=at position .5 with {\draw (0,-0.1) -- (0,0.1);}},line width=2.000]

\DeclareFieldFormat{postnote}{#1}
\DeclareFieldFormat{multipostnote}{#1}
\DeclareFieldFormat{pages}{#1}
\title{Root group data (RGD) systems of affine type for significant subgroups of isotropic reductive groups over $k[t,t^{-1}]$} 
\author{Yuan Zhang}
\date{}
\pgfplotsset{compat=1.18}

\begin{document}
    \maketitle
    \begin{abstract}
        Given a connected isotropic reductive not necessarily split $k$-group $\mathcal{G}$ with irreducible relative root system, we construct root group data (RGD) system of affine type for significant subgroups of $\mathcal{G}(k[t,t^{-1}])$, which can be extended to the whole group $\mathcal{G}(k[t,t^{-1}])$ under certain additional requirements. We rely on the relative pinning maps from \cite{cite:ele_subgp_of_iso_red_gp} to construct the affine root groups. To verify the RGD axioms, we utilize the properties of the affine root groups, and the properties of reflections associated with the $k$-roots of $\mathcal{G}$. 
    \end{abstract}
    \section{Introduction}
        \input{sections/Introduction.tex}

    \section{Review on root group data (RGD) system and affine roots}\label{sec:Review on root group data (RGD) system and affine roots}
        \input{sections/Review_on_RGD_and_affine_roots.tex}
    \section{Review on reductive groups}\label{sec:Review on reductive groups}
        \input{sections/Review_on_reductive_groups.tex}
    \section{Summary of needed results}\label{sec:Summary of needed results}
        \input{sections/Summary_of_needed_results.tex}

    \section{Main theorem and its proof}\label{sec: Main theorem and its proof}
        \input{sections/Main_theorem_and_its_proof.tex}
    \section{Examples}\label{sec:Examples}
        \input{sections/Examples.tex}
    \section{Acknowledgement}
        This paper consists of results during my PhD studies at University of Virginia, and, so, I thank my PhD advisor Peter Abramenko's continuous support. 
    \printbibliography
\end{document}

%% file: sections/Introduction.tex
\subsection{Motivation and history}

BN pairs were first introduced by Jacques Tits in the study of certain matrix groups (in particular, reductive and semisimple groups). Through introducing buildings (first in \cite{cite:building_first_1955}), Tits generalized the projective geometric space, and he utilized BN pairs to reach a unified simplicity proof when studying Chevalley groups over fields. 

It is then that the link between BN pairs and buildings was established, introducing discrete geometric ideas into abstract group theory. BN pairs arise from root groups, such as in the Chevalley groups over a field $k$, the root groups are the 1-parameter subgroups $U_{\alpha}$ (isomorphic to the additive group of $k$ as groups) where $\alpha$ is a root in the corresponding root system of the Chevalley group. To make the connection between root groups and BN pairs, Tits formulated the notion of a \textbf{root group data (RGD) system} that occurs in algebraic groups. This system always gives rise to a BN pair. 

With every BN pair, there is an associated Weyl group (and Coxeter system). In the case of the Chevalley groups (and, more generally, the isotropic reductive groups) over fields, the Weyl groups are finite, and give rise to spherical buildings. 
Then, Iwahori and Matsumoto (in \cite{cite:Iwa_Matsu_p-adic_Ch_gp}) used Chevalley groups over a local field to introduce BN pairs (hence buildings) with infinite Weyl groups $W_{\text{aff}}(\Psi)$ of affine type (where $\Psi$ is the root system corresponding to the group). In \cite{cite:gp_rd_sur_corp_local}, Tits and Bruhat then generalized this work and provided affine BN pairs for isotropic reductive groups over local fields, which again arose from RGD systems. As an example, $k(t)$ has discrete valuations $v^{+}$ (with prime element $t$) and $v^{-}$ (with prime element $t^{-1}$), this fact gives us two different affine BN pairs on $\mathcal{G}(k(t))$ for a Chevalley group $\mathcal{G}$. Since $k[t,t^{-1}]$ is dense in $k(t)$ (in the topology induced by the two discrete valuations), the two affine BN pairs on $\mathcal{G}(k(t))$ further induce two affine BN pairs on $\mathcal{G}(k[t,t^{-1}])$. These two affine BN pairs on $\mathcal{G}(k[t,t^{-1}])$ give rise to a twin building on the group. 

Given a Chevalley group $\mathcal{G}$, $G=\mathcal{G}(k[t,t^{-1}])$ can be seen as the $k[t,t^{-1}]$-points of a linear algebraic group; and as a Kac-Moody group over $k$ of affine type (in the sense of Tits'), as well. In \cite{cite:Uniq_and_pres_Kac_gp_over_field}, Tits provided an RGD system of not necessarily spherical type for Kac-Moody groups of ``split'' type over fields (in the sense of Tits', which covers the case of $G$ we just described), which system give rise to twin BN pairs and, hence, twin buildings. This result has many applications, chief among which is the study of action on twin buildings. 

However, for an isotropic reductive $k$-group $\mathcal{G}$, $G=\mathcal{G}(k[t,t^{-1}])$ is not necessarily known to be a Kac-Moody group over $k$ in the sense of Tits'. 
In \cite[sec: 3.2]{cite:tw_bd_kac}, Tits provided a twin BN pair in this more general situation when $\mathcal{G}$ is almost simple and simply connected, however, an explicit construction of an RGD system for this case is not present in the literature to the best of my knowledge\footnote{It later came to my attention that there is a non-explicit construction of RGD system for $k[t,t^{-1}]$-points of reductive $\mathbb{F}_{q}$-groups for fields $k$ that contain $\mathbb{F}_{q}$ through the method of Galois descending, see \cite[prop: 10.2, 10.3]{cite:fini_prop_Ch_poly}.}. 
This fact brings us to the \textbf{main goal} of this paper: to provide an explicit construction of an RGD system in the not necessarily split reductive case for some significant subgroups. 

\subsection{Summary}

We start with a review on RGD systems and affine roots in Section \ref{sec:Review on root group data (RGD) system and affine roots}. Then we review the reductive groups in Section \ref{sec:Review on reductive groups}. In Section \ref{sec:Summary of needed results}, we summarize the construction and relevant properties of \textbf{relative pinning maps} introduced by Petrov and Stavrova in \cite{cite:ele_subgp_of_iso_red_gp}. 
In Section \ref{sec: Main theorem and its proof}, 
we first utilize the relative pinning maps to construct the \textbf{affine root group} $U_{\alpha_{a',l}}$ associated with an affine root $\alpha_{a',l}\in \Phi$, where $\Phi$ is the set of affine roots associated with the relative root system of isotropic reductive $k$-group $\mathcal{G}$ (see \ref{n:Construction of Affine root group in non-split case} for construction of affine root groups). 
Then, with properties of relative pinning maps (see \ref{n:Properties of realtive pinning maps}), properties of absolute pinning isomorphisms (see \ref{n:Existence of needed faithful rational representation for the split case}), and properties of reflections in the co-space in respect to the unique dual of a k-root (see \ref{d:k-root of reductive group}), we prove the main theorem (Theorem \ref{t:ref{t:Construction of RGD system for the non-split case} without further requirement}): \bigskip

\begin{thm*}
    \label{n:Summary of main theorem}
    For a connected isotropic reductive affine algebraic group $\mathcal{G}$ with irreducible relative root system, where we denote by $\mathcal{G}(k[t,t^{-1}])^{+}$ the elementary subgroup (see definition of elementary subgroup at \ref{n:Important proteries of generalize unipotent elements and affine root group}(d)), we have that 
    \[\Bigl(\mathcal{G}\bigl(k[t,t^{-1}]\bigr)^{+}C_{\mathcal{G}}(S)(k),(U_{\alpha_{a',l}})_{\alpha_{a',l}\in \Phi},C_{\mathcal{G}}(S)(k)\Bigr)\] 
    is an RGD system of affine type. 
\end{thm*}
Furthermore, we can extend the main theorem (Remark \ref{rmk:extending main theorem}): \bigskip 
\begin{remark*}
    \label{n:Summary of remark ref{rmk:extending main theorem}}
    It is a consequence of \cite[cor:6.2]{cite:hom_inv_of_non_stable_K1} that
    the main theorem can be extended to $\mathcal{G}(k[t,t^{-1}])$ from $\mathcal{G}(k[t,t^{-1}])^{+}C_{\mathcal{G}}(S)(k)$ under additional hypotheses. 
\end{remark*}

In Section \ref{sec:Examples}, we utilize the main theorem to construct RGD systems of affine type for Larurent polynomial points of split reductive $k$-groups (see \ref{ex:Split case}), and for special unitary groups (see \ref{ex:Unitary group case}).

%% file: sections/Review_on_RGD_and_affine_roots.tex
\subsection{Root group data (RGD) systems}

To review the \textbf{roots of Coxeter complex}, one refers to \cite[sec:3.4]{cite:bd}, \cite[5.5.4]{cite:bd}, and \cite[sec:5.8.5]{cite:bd}. \bigskip 

\begin{defn}[Prenilpotent pair of roots]\label{d:Prenilpotent pair of roots}
    (\cite[8.41]{cite:bd}): Consider the Coxeter system (W,S). Let $\Sigma$ be the simplicial Coxeter complex of type (W,S). Let $\Phi$ be the set of roots of $\Sigma$. 
    Given $\alpha,\beta\in \Phi$, the pair $\{\alpha,\beta\}$ is called \textbf{prenilpotent} if $\alpha\cap \beta$ and $(-\alpha)\cap(-\beta)$ each contain at least one chamber. In this case, we set 
    \[[\alpha,\beta]:=\Bigl\{\gamma\in\Phi | \gamma\supset \alpha\cap\beta\;and\;-\gamma\supset(-\alpha)\cap(-\beta)\Bigr\}.\] 
    The left and(or) right ``$[]$'' is replaced by ``$()$'' by excluding $\alpha$ and(or) $\beta$ from the set $[\alpha,\beta]$. 
\end{defn}

Some results regarding prenilpotent pairs of roots can be found in \cite[ch:8.5.3]{cite:bd}. Specifically, \cite[8.44]{cite:bd} helps in understanding in affine roots, and \cite[8.42]{cite:bd} is helpful. 

For ease of reference, we state the axioms for a general RGD system, we follow closely to \cite[sec: 8.6.1]{cite:bd} for this purpose: 
\begin{defn}[General RGD system]\label{d:General RGD system}
    Let $(W,S)$ be an arbitrary Coxeter system. Let $\Sigma:=\Sigma(W,S)$ be its Coxeter complex. We identify $\mathcal{C}(\Sigma)$ (i.e. the set of chambers of $\Sigma$) with $W$ (see \cite[5.66]{cite:bd} for this identification) and take $1$ to be the fundamental chamber. Let $\Phi$ denote the set of roots
    of $\Sigma$ (we denote by $\Phi_{+}$ (resp. $\Phi_{-}$) the positive (resp. negative) set of roots defined by containing (resp. not containing) 1). 

    Recall that, for $s\in S$, $\alpha_{s}$ denotes the roots characterized by 
    \[\mathcal{C}(\alpha_{s})=\Bigl\{w\in W \Big| l(sw)>l(w)\Bigr\}.\] 
    We write $U_{s}$ to denote $U_{\alpha_{s}}$, and $U_{-s}$ to denote $U_{-\alpha_{s}}$. 

An \textbf{RGD system of type $(W,S)$} is a triple $(G,(U_\alpha)_{\alpha\in\Phi},T)$, where $T$ and the $U_\alpha$'s are
subgroups of $G$, satisfying the following set of axioms: 

\begin{enumerate}
    \item[(RGD0)] For all $\alpha\in\Phi$, $U_\alpha\neq 1$. 
    \item[(RGD1)] For all $\alpha\neq \beta$ in $\Phi$ such that $\{\alpha,\beta\}$ is prenilpotent, we have \[[U_{\alpha}, U_{\beta}]\leq U_{(\alpha,\beta)}.\]  
    \item[(RGD2)] For all $s\in S$, there is a function $m:U_s^*\to G$ such that, for all $u,v\in U_s^*$, and all $\alpha\in\Phi$, we have 
    \[m(u) \in U_{-s} u U_{-s} \text{, } m(u)U_\alpha m(u)^{-1}=U_{s\alpha}\;,\text{ and }\] 
    \[m(u)^{-1}m(v)\in T.\] 
    \item[(RGD3)] For all $s\in S$, $U_{-s}\not\leq U_+:=\langle U_\alpha | \alpha\in\Phi_+\rangle $. (where $U_{-}$ denotes $\langle U_{\alpha}|\alpha\in\Phi_{-} \rangle$). 
    \item[(RGD4)] We have that $G=T\langle U_\alpha | \alpha\in\Phi\rangle $.
    \item[(RGD5)] The subgroup $T$ normalizes each $U_\alpha$, i.e. 
    $\displaystyle T\subset \bigcap_{\alpha\in \Phi}N_{G}(U_{\alpha})$.  
\end{enumerate}

\textbf{Note}: We can replace (RGD3) with either of the following items and receive an equivalent set of axioms: 
\begin{enumerate}
    \item[(RGD3'')] For all $s\in S$, we have that $U_{-s}\not\subset U_{+}\;\&\; U_{s}\not\subset U_{-}$.  
    \item[(RGD3)*] $TU_{+}\cap U_{-}=\{1\}$.   
\end{enumerate}
See \cite[8.77]{cite:bd} for (RGD3'') and ``hint 2'' on \cite[p.404]{cite:gp_act_tb} for (RGD3)*. 

When we speak of an \textbf{RGD system of affine type}, we mean that the Weyl group $W$ is an affine Weyl group. See review on affine roots in the following. 
\end{defn}\bigskip

\subsection{Affine roots}
    
We gather some background on affine roots. A good reference is \cite{cite:ref_gp_and_cox_gp}. Beyond just recording the foundations, we bring the foundations in \cite{cite:ref_gp_and_cox_gp} and the results in \cite{cite:bd} together to obtain some tools we need. 
\bigskip

\begin{note}[Affine roots and affine Weyl group]\label{n:Affine roots and affine Weyl group}
    Let $\Psi$ be a 
    \textbf{root system} in Euclidean space $V:=\mathbb{R}^{n}$. Let $\Pi=\{a_{1},\cdots,a_{n}\}$ be a base of $\Psi$ (see \cite[ch:III]{cite:intro_lie_rep} for a review on root systems). We use $(\cdot,\cdot)$ to denote 
    the standard inner product associated with the Euclidean space. We denote by $\Psi_{\pm}$ the sets of $\pm$ roots, respectively, and by $a_{0}$ a \textbf{root of maximal height} in $\Psi_{+}$ (having a unique highest root requires the \textbf{irreducible} condition on $\Psi$). 
    In this setup, we construct the \textbf{set of affine roots associated with $\Psi$} to be 
    \[\Phi=\Bigl\{\alpha_{a,l}|a\in\Psi,\;l\in\mathbb{Z}\Bigr\},\text{ where }\alpha_{a,l}:=\Bigl\{v\in V|(a,v)\geq -l\Bigr\}.\]  
    An observation is that, by this construction, $\alpha_{a,l}=-\alpha_{-a,-l}$. 
    It can be seen that there are \textbf{affine hyperplanes} (associated with affine roots) $\partial\alpha_{a,l}=\{v\in V|(a,v)=-l\}$. We define the reflection $s_{a,l}$ to be the \textbf{affine reflection} with respect to the affine hyperplane $\partial\alpha_{a,l}$. 
    We also write $H_{a,r}:=\{v\in V|(a,v)=-r\}$ for $a\in\Psi$ and $r\in \mathbb{R}$. 
    Note that $H_{a,r}=\partial\alpha_{a,r}$ if $r\in \mathbb{Z}$. 
    We denote by $W_{\text{aff}}(\Psi):=\langle s_{a,l}|a\in\Psi, l\in\mathbb{Z} \rangle $ the \textbf{affine Weyl group} of $\Psi$. 
    We define the \textbf{simple affine roots} $\alpha_{0}:=\alpha_{-a_{0},1}$ and $\alpha_{i}:=\alpha_{a_{i},0}$ for all $1\leq i\leq n$. 
    We also write $s_{i}$ for the affine reflection with respect to $\alpha_{i}$ for $0\leq i\leq n$ respectively. 
    By convention, we define $<b,a>:=\displaystyle\frac{2(a,b)}{(a,a)}$ for $a,b\in \Psi$, with the \textbf{dual root} of $a\in \Psi$ being $a^{\vee}:=\displaystyle\frac{2a}{(a,a)}$.  
    By the requirement of reduced root system that the only possible scalar multiples of a root $a$ are $\pm a$, 
    determining $e_{a}:=\frac{a}{\sqrt{(a,a)}}$ is equivalent to determining $a\in \Psi$ in a reduced root system. 
    Therefore, we have a one-to-one correspondence (whenever $H_{a,r}$ exists) 
    of  
    \[(a,r)\leftrightarrow H_{a,r}\] 
    for a root $a\in\Psi$ and for $r\in \mathbb{R}$.  
\end{note}\bigskip

\begin{note}\label{n:A Summary of results about affine roots}
    We give a \textbf{summary of results regarding affine roots and affine Weyl groups}. We list the specific context in front of the statement. If there is no specific context, then the statement apply in the most general case: 
    \begin{enumerate}
        
        \item The reflection $s_{a,l}$ fixes $\partial\alpha_{a,l}$ pointwise and sends 0 to $-la^{\vee}$. Hence, we may write $s_{a,l}= \tau(-la^{\vee})\circ s_{a,0}$, where $\tau(\lambda)$ is the translation $\tau(\lambda):v\mapsto v+\lambda$ for $v,\lambda\in \mathbb{R}^{n}$. (\cite[sec:4.1]{cite:ref_gp_and_cox_gp}). 
        
        \item It can be seen with routine calculations that   
        \[s_{a,l}(\alpha_{b,m})=\alpha_{s_{a,0}(b),m-l<b,a>}.\]  
        
        \item (\textbf{Requires the irreducible condition}) 
        We consider the positive affine roots to be the ones that contain the fundamental chamber. Then, we have  
        \[\Phi_{+}=\Bigl\{\alpha_{a,l}\in \Phi \Big| (a\in\Psi_{+}\;and\; l\geq 0) \text{ or } (a\in\Psi_{-}\;and\;l\geq 1)\Bigr\}\;,\text{ and }\]
        \[\Phi_{-}=\Bigl\{\alpha_{a,l}\in\Phi \Big| (a\in\Psi_{+}\&l\leq -1)\text{ or }(a\in\Psi_{-}\&l\leq 0)\Bigr\}.\] 
        
        \item \textbf{(Requires the reduced condition)} The pair $\{\alpha_{a,l},\alpha_{b,m}\}\subset \Phi$ is prenilpotent if and only if $e_{a}\neq -e_{b}$ (see \cite[8.44]{cite:bd}). 
        Then, according to the reduced condition on a root system, it can be seen that  
        \[\text{the pair $\{\alpha_{a,l},\alpha_{b,m}\}\subset \Phi$ is prenilpotent if and only if $a\neq -b$}.\] 
        
        \item[4'.] \textbf{(After relaxing the reduced condition of item 4)} When the reduced condition on $\Psi$ is relaxed, we instead have that $e_{a}\neq -e_{b}$ if and only if $ka\neq -nb$ for all $k,n\in \mathbb{Z}_{>0}$ (note that, by the construction of root systems, the only possible sets of proportional roots are $\{\pm a\}$, $\{\pm a,\pm \frac{1}{2}a\}$, and $\{\pm a, \pm 2a\}$; see \cite[14.7]{cite:borel}). Therefore,  
        \[\text{$\alpha_{a,l}$ and $\alpha_{b,m}$ are prenilpotent if and only if $ka\neq -nb$ for any $k,n\in \mathbb{Z}_{>0}$}.\] 
        
        \item 
        When $\{\alpha_{a,l},\alpha_{b,m}\}\subset \Phi$ is prenilpotent: 
        \[\Bigl\{\alpha_{pa+qb,pl+qm}\in\Phi \Big| p,q\in\mathbb{R}_{\geq 0}\Bigr\}\subset [\alpha_{a,l},\alpha_{b,m}].\]  
        To see that $\{\alpha_{pa+qb,pl+qm}\in\Phi|p,q\in\mathbb{R}_{\geq 0}\}\subset [\alpha_{a,l},\alpha_{b,m}]$, note that, for $v\in V$ and $p,q\in\mathbb{R}_{\geq 0}$: 
        \[(a,v)\geq -l\text{ and }(b,v)\geq -m \text{ implies } p(a,v)+q(b,v)\geq -pl-qm\;,\text{ and }\]
        \[ (a,v)\leq -l\text{ and }(b,v)\leq -m \text{ implies } p(a,v)+q(b,v)\leq -pl-qm.\] 
        Thus, for $\alpha_{pa+qb,pl+qm}\in\{\alpha_{pa+qb,pl+qm}\in\Phi|p,q\in\mathbb{R}_{\geq 0}\}$, we have $\alpha_{a,l}\cap\alpha_{b,m}\subset \alpha_{pa+qb,pl+qm}$ and $-\alpha_{a,l}\cap-\alpha_{b,m}\subset -\alpha_{pa+qb,pl+qm}$. 

        Under the further assumption that $\Psi$ be crystallographic and reduced, we can see that $[\alpha_{a,l},\alpha_{b,m}]=\{\alpha_{pa+qb,pl+qm}\in\Phi|p,q\in\mathbb{Q}_{\geq 0}\}$. To do the converse inclusion, we need \cite[8.45]{cite:bd}, and we develop into the two cases according to \cite[8.45]{cite:bd}'s two different cases to do so. 
        We do not utilize this extended result, so the proof is not included.

    \end{enumerate}
\end{note}

%% file: sections/Review_on_reductive_groups.tex
We opt for the classical setting of algebraic geometry in this paper. There are good references for algebraically closed fields. However, since we are considering relative cases, we need to understand things beyond algebraically closed fields. Some good references for this purpose are \cite[sec:1.1]{cite:Reducktive_Gruppen_AP} and \cite[Ch: XII]{cite:LAG_hum}. 

We provide a summary on algebraic groups, and one on reductive groups. All the following results and concepts are known. The main reference is \cite{cite:borel}, with alternative references \cite{cite:LAG_hum}, \cite{cite:TAS_LAG}, and \cite{cite:LAG_PA}.

We refer to \cite{cite:borel} for information on \textbf{algebraic groups} (we only consider affine algebraic groups). 
The statement \cite[1.10]{cite:borel} tells us we always have a \textbf{faithful $k$-rational representation} of $k$-group $\mathcal{G}$ (see the definition of a \textbf{rational representation} in \cite[1.6(8)]{cite:borel}). 
We define \textbf{tori} as in \cite[8.5]{cite:borel}, and \textbf{split tori} as in \cite[8.2]{cite:borel} (or \cite[Definition 5]{cite:Reducktive_Gruppen_AP}). In particular, note that a connected $k$-group $\mathcal{G}$ must contain a maximal torus that is defined over $k$ (see \cite[18.2]{cite:borel}). Then, we may choose a maximal torus of a connected $k$-group to be defined over $k$. 
A review of the construction of the \textbf{Lie algebra of an algebraic group $\mathcal{G}$}, denoted by $\mathcal{L}(\mathcal{G})$, (according to derivation) can be found at \cite[3.3]{cite:borel}. Alternatively, one can construct $\mathcal{L}(\mathcal{G})$ as $T(\mathcal{G})_{e}$ (tangent space at the identity element), as described in \cite[3.2.1]{cite:LAG_PA} from the classical point of view. 

Characters and cocharacters are important tools in this paper, so we discuss them below: 
\bigskip

\begin{defn}[Character]\label{d:Characters}
    (\cite[5.2]{cite:borel}): Let $\mathcal{G}$ and $\mathcal{G}'$ be $k$-groups. We write $\text{Mor}(\mathcal{G},\mathcal{G}')$ for the algebraic group morphisms on $\mathcal{G}$ to $\mathcal{G}$, and $\text{Mor}(\mathcal{G},\mathcal{G}')_{k}$ for the set of such morphisms that are defined over $k$. 
    We write $X(\mathcal{G}):=\text{Mor}(\mathcal{G},\mathbf{GL_{1}})$ and call its elements \textbf{characters} of G. 
    Since $\mathbf{GL}_{1}$ is commutative, $X(\mathcal{G})$ becomes an abelian group by $(a_1+a_2)(g)=a_1(g)a_2(g)$ for $a_{1},a_{2}\in X(\mathcal{G})$. 
    
\end{defn}\bigskip
\begin{defn}[Cocharacter]\label{d:The multiplicative one parameter subgroups}
    (\cite[8.6]{cite:borel}): For a $k$-group $\mathcal{G}$, consider \[X_{*}(\mathcal{G})=Mor(\mathbf{GL_{1}},\mathcal{G}).\]  
    The elements of $X_{*}(\mathcal{G})$ are called \textbf{cocharacters}. 
    We have a map \[X(\mathcal{G})\times X_{*}(\mathcal{G})\to \mathbb{Z}=X(\mathbf{GL_{1}})\] 
    defined by $<\gamma,\chi>:=m$ if $(\gamma\circ\chi)(x)=x^{m}$ for $x\neq 0$. 
    Then, $X_{*}(\mathcal{G})$ (akin to $X(\mathcal{G})$) becomes an abelian group if $\mathcal{G}$ is commutative: we set $(\chi_{1}+\chi_{2})(\lambda)=\chi_{1}(\lambda)\chi_{2}(\lambda)$ for $\chi_{1},\chi_{2}\in X_{*}(\mathcal{G})$. 
    Hence, the above map becomes a bilinear map of abelian groups when $\mathcal{G}$ is commutative. 
    For a torus T, we have  
    \[X(T)\times X_{*}(T)\to \mathbb{Z}\;,\] 
    which we often denote by $<\cdot,\cdot>$. For more information, consult \cite[sec: 3.2]{cite:TAS_LAG}.
\end{defn}

As in \cite[8.17]{cite:borel}, we construct the \textbf{roots} of algebraic group $\mathcal{G}$ with respect to an arbitrary torus T according to the weight spaces $\mathfrak{g}_{\alpha}$'s of the adjoint representation (see about the adjoint representation in \cite[3.13]{cite:borel} or \cite[3.4.1]{cite:LAG_PA}). We denote the \textbf{set of roots} of $\mathcal{G}$ \textbf{with respect to} T by $\Phi(\mathcal{G},T)$. \bigskip

\begin{remark}[Lie algebras of Chevalley group]\label{rmk:Lie algebra of Chevalley group}
    According to \cite[Ch:2, Corollary 1]{cite:st}, \cite[Ch:2, Corollary 3]{cite:st}, and the remark above \cite[Ch:5, Corollary 1]{cite:st}, we see that, for an algebraically closed field $K$, the Lie algebra of the Chevalley group $\mathcal{G}(K)$ is $\mathcal{L}_{\mathbb{Z}}\otimes_{\mathbb{Z}}K$, where $\mathcal{L}_{\mathbb{Z}}$ is defined in \cite[Ch:2, Corollary 2]{cite:st}. 
\end{remark}

We define \textbf{reductive groups} as in \cite[11.21]{cite:borel}. We discuss some important properties of reductive groups. 

\begin{defn}[Split reductive group and admissible isomorphism (i.e. absolute pinning isomorphism)]\label{d:Split reductive group and admissible isomorphism}
    (\cite[18.6-18.7]{cite:borel} and \cite[8.1.1]{cite:TAS_LAG}): 
    Let $\mathcal{G}$ be a reductive and connected $k$-group. Let T be a maximal torus of $\mathcal{G}$. Write $\Phi:=\Phi(\mathcal{G},T)$. Let $K$ be an algebraically closed field containing $k$. 
    For each $\alpha\in \Phi$, there is a unique connected subgroup $U_{\alpha}$ of G that is normalized by T and has the property that $\mathcal{L}(U_{\alpha})=\mathfrak{g}_{\alpha}$. We have an isomorphism $x_{\alpha}:Add(K)\xrightarrow{\cong} U_{\alpha}$ such that 
    \[tx_{\alpha}(\lambda)t^{-1}=x_{\alpha}\Bigl(\alpha(t)\lambda\Bigr)\] 
    for all $t\in T$, and for all $\lambda\in K$. It can be seen that $U_{\alpha}$ has dimension 1 (see more on this in \cite[13.18]{cite:borel}). 
    We call these $x_{\alpha}(\cdot)$'s \textbf{admissible isomorphisms} (as in\cite{cite:LAG_hum}), or as \textbf{absolute pinning isomorphisms} (As \cite{cite:ele_subgp_of_iso_red_gp} refers them as a part of ``\'epinglage'' originated in \cite[Exp XXIII, def 1.1]{cite:sche_en_gps}, while we use term ``absolute'' to emphasize the fact that $a$ is an absolute root). We call the $U_{\alpha}$'s \textbf{absolute root groups}.

    We say that $\mathcal{G}$ is \textbf{split over} $k$ (or just \textbf{split} when $k$ is clear by context) if we can choose a torus $T$ to be split over $k$, as well as isomorphisms $x_{\alpha}$ to be defined over $k$. Furthermore, it can be checked that a reductive group $\mathcal{G}$ is split over $k$ if it has a maximal torus that splits over $k$ (see \cite[18.7]{cite:borel}). Then, by \cite[8.11]{cite:borel}, we can see that $\mathcal{G}$ is always split over some finite separable extension $\hat{k}$ of $k$. 
    By the equation in the proof of \cite[7.3.3]{cite:TAS_LAG}, it can be seen that each $\mathcal{U}_{\alpha}\subset [\mathcal{G},\mathcal{G}]$. Furthermore, note that $[\mathcal{G},\mathcal{G}]$ is semisimple and $\mathcal{G}=[\mathcal{G},\mathcal{G}]Z(\mathcal{G})^{\circ}$ (see \cite[8.1.6]{cite:TAS_LAG} or \cite[14.2]{cite:borel}). This result makes the Chevalley group relevant to our study, as the Chevalley groups are precisely the split semisimple groups. 
    It can also be seen (as in \cite[13.17]{cite:borel}) that $C_{\mathcal{G}}(T)=T$. 

\end{defn}\bigskip

\begin{defn}[Reduced root system in a reductive group]\label{n:Reduced root system in reductive group}
    (\cite[14.8]{cite:borel}): For a connected reductive group $\mathcal{G}$, and a maximal torus T of $\mathcal{G}$, $\Phi:=\Phi(\mathcal{G},T)$ is a reduced root system which we call the \textbf{absolute root system}. 

    Note that, for $\alpha\in \Phi$, we take $\alpha^{\vee}\in X_{*}(T)\subset V^{\vee}:=\mathbb{R}\otimes_{\mathbb{Z}}X_{*}(T)$ (see \cite[7.3.5]{cite:TAS_LAG} for justification for the element containment) to be the \textbf{unique dual}, as in \cite[7.1.8]{cite:TAS_LAG} (we also call $\alpha^{\vee}$ the \textbf{coroot} of $\alpha$, since $\alpha^{\vee}$ is a cocharacter). The reflection in $X(T)\otimes_{\mathbb{Z}}\mathbb{R}$ can be written as \[s_{\alpha}(v)=v-<v,\alpha^{\vee}>\alpha\]

    for $v\in V:=\mathbb{R}\otimes_{\mathbb{Z}} X(T)$ and $\alpha\in \Phi\subset V$. Note that the pairing $<\cdot,\cdot>$ of $X(T)$ and $X_{*}(T)$ induces a pairing of $V$ and $V^{\vee}$, which we also denote by $<\cdot,\cdot>$. 

    We note that there are also reflections in $V^{\vee}$, which we refer to as \textbf{reflections in the co-space}, given by 
    \[s_{\alpha^{\vee}}(v^{*})=v^{*}-<\alpha,v^{*}>\alpha^{\vee}.\]  
    The reflection in the co-space gives $s_{\alpha}(\beta)^{\vee}=s_{\alpha^{\vee}}(\beta^{\vee})$ for $\alpha,\beta\in \Phi$ (with the bilinearity of $<\cdot,\cdot>$). In fact, there are stronger properties satisfied by the quadruple $(\Phi,X(T),\Phi^{\vee},X_{*}(T))$: it is a root datum (see \cite[7.4.3]{cite:TAS_LAG} and \cite[7.4.1]{cite:TAS_LAG}, where more information about reflections in the co-space is also mentioned).

    We can consider the cocharacter group $X_{*}(T)$ identified to a lattice $\{\lambda\in V|(\lambda,X(T))\subset \mathbb{Z}\}$ in $V$ with correspondence the $\alpha^{\vee}\leftrightarrow \frac{2\alpha}{(\alpha,\alpha)}$ for root $\alpha$ (see \cite[15.3.6]{cite:TAS_LAG}, the bilinear form $(\cdot,\cdot)$ here is constructed in \cite[7.1.7]{cite:TAS_LAG}). 
    Then we may define, in this context, that $<\beta,\alpha>:=<\beta,\alpha^{\vee}>$ with equality $<\beta,\alpha>=(\beta,\frac{2\alpha}{(\alpha,\alpha)})$ (to be consistent with the conventional notation when studying root systems). We can do so because each $\alpha$ has a unique $\alpha^{\vee}$, and because of how $\alpha^{\vee}$ is taken in \cite[7.1.8]{cite:TAS_LAG}. We use both notations interchangeably. Note that, $\alpha\mapsto \alpha^{\vee}$ is not (necessarily) linear, so $<\cdot,\cdot>$, when both coordinates are roots, need not be bilinear.

\end{defn}\bigskip

\begin{defn}[$k$-root of reductive group]\label{d:k-root of reductive group}
    (\cite[21.1]{cite:borel}): 
    Let $S$ be a maximal $k$-split torus of a connected reductive $k$-group $\mathcal{G}$. 
    Denoted by $_{k}\Phi:=\Phi(\mathcal{G},S)$ the \textbf{set of $k$-roots of $\mathcal{G}$ with respect to $S$}. 
    It can be seen (\cite[21.6-21.7]{cite:borel}) that $_{k}\Phi$ is a not necessarily reduced root system. Hence, we call it a \textbf{relative root system}, or a \textbf{$k$-root system}. More details about the relative root system can be found in \cite[sec: 15.3]{cite:TAS_LAG}. In particular, the notions of $<\cdot,\cdot>$ and coroots generalize to the context of $k$-roots, and the reflections and reflections in the co-space are described in \cite[15.3.8]{cite:TAS_LAG} (the generalization of the formulae in relative case ``are the same'' as those in \ref{n:Reduced root system in reductive group}). 
\end{defn}\bigskip 

\begin{note}[Existence of needed faithful rational representation for the split case]\label{n:Existence of needed faithful rational representation for the split case}
    Consider a connected reductive split $k$-group $\mathcal{G}$ with a maximal torus $T$ that is $k$-split. 
    Write $\Phi=\Phi(\mathcal{G},T)$. 
    Then, $\mathcal{G}$ shares the absolute root system, and hence the absolute root groups, with its commutator subgroup $[\mathcal{G},\mathcal{G}]$ (which is semisimple). For a fixed absolute root $a$, and for any two absolute pinning isomorphisms $x_{a}(\lambda)$ and $y_{a}(\lambda)$, 
    we must have that $x_{a}(\lambda)=y_{a}(c\lambda)$ for some $c\in \overline{k}^*$ (see \cite[18.6]{cite:borel}). We then see that the $x_{\alpha}(\cdot)$'s Steinberg defined in \cite[Ch:3]{cite:st} satisfy the requirement needed for absolute pinning isomorphisms in the reductive case. 
    As a consequence, absolute pinning isomorphisms inherit the properties of construct by Steinberg, in particular: 

    There exists a faithful rational representation defined over $k$, $\rho:G\to GL_{N}(K)$, such that the elements of $\rho(\mathcal{U}_{\alpha})$ are uni-upper-triangular (resp. uni-lower-triangular) for $\alpha$ positive root (resp. negative root) in $\Phi$. 
\end{note}

Our main concern is the not necessarily split case. \bigskip

\begin{note}[Set up for the not necessarily split case]\label{n:Set up for non-split case} 
    The notations and setup follow \cite[sec: 1.2.3]{cite:Reducktive_Gruppen_AP} and \cite[sec: 21]{cite:borel}. 

    \begin{enumerate}
        \item Let $\mathcal{G}$ be a connected reductive $k$-group that does not split. 
        Consider a maximal $k$-split torus $S$ of $\mathcal{G}$; we can always find a maximal tours $T$ of $\mathcal{G}$ that is defined over $k$ and contains $S$ (see \cite[Satz 13]{cite:Reducktive_Gruppen_AP}). Consider the $k$-roots in $_{k}\Phi:=\Phi(\mathcal{G},S)$ and the roots in $\Phi:=\Phi(\mathcal{G},T)$. The inclusion of $S\hookrightarrow T$ induces a group projection $j:X(T)\rightarrow X(S)$ by restriction $\chi\mapsto \chi|_{S}$. Then, for $a,b\in \Psi:=\Phi(\mathcal{G},T)$, we have that $<b,j(a)>=<j(b),j(a)>$, as $j(b):=b|_S$ 
        
        As shown in \cite[21.8]{cite:borel},  
        \[\text{$_{k}\Phi\subset j(\Phi)\subset\; _{k}\Phi\cup \{0\}$ and $_{k}\Phi_{+}\subset j(\Phi_{+})\subset\; _{k}\Phi_{+}\cup \{0\}$}.\] 
        We also write 
        \[\eta(a'):=j^{-1}(a')\cap \Phi\] 
        for $a'\in\; _{k}\Phi$. 

        \item 
    
        For $a'\in\; _{k}\Phi$, we write \[(a'):=\{\lambda a'| \lambda\in \mathbb{Z}_{>0}\}.\] 
        Then, we have that \[(a')=
            \{a'\} \text{ or } \{a', 2a'\} 
        .\]  

        \item For $a'\in\; _{k}\Phi$, we construct  
        \[U_{a'}:=U_{(a')}:=\Bigl\langle U_{a}|a\in \eta\bigl((a')\bigr) \Bigr\rangle \;,\] 
        where $U_{a}$ for $a\in \Phi$ are constructed in \ref{d:Split reductive group and admissible isomorphism}. We can observe that $(a')$ is ``special'' in $_{k}\Phi$ (``special'' in the sense of \cite[14.5]{cite:borel}: see more at \cite[14.7]{cite:borel} and \cite[21.7]{cite:borel}). We infer that $\eta((a'))$ is ``special'' in $\Phi$.  
        By \cite[14.5]{cite:borel}, we have that $U_{a'}$ is \textbf{directly spanned} (see \cite[14.3]{cite:borel}) in any order by $(U_{a})_{a\in \eta((a'))}$, that is, for any ordering of $\{a_{1},\cdots,a_{l}\}=\eta((a'))$, the map  
        \[U_{a_{1}}\times \cdots \times U_{a_{l}}\to U_{a'} \text{ by } (u_{1},\cdots,u_{l})\mapsto u_{1}\cdots u_{l}\] 
        is an isomorphism of affine varieties. 

        \item (\cite[21.10]{cite:borel}): Let $a'\in\; _{k}\Phi$. A key observation is that, if $(a')=\{a'\}$, then $U_{a'}$ is commutative. If $(a')=\{a',2a'\}$, then $U_{a'}$ is no longer commutative. However, we have $U_{2a'}=Z(U_{a'})$ (see \cite[4.10]{cite:hom_ag_sim}). 
    \end{enumerate}
\end{note}

%% file: sections/Summary_of_needed_results.tex
\begin{note}[Context]\label{n:Context in needed results}

    Let $S$ be a maximal $k$-split torus in a connected isotropic reductive $k$-group $\mathcal{G}$. Let $T$ be a maximal torus that contains $S$, and be defined over $k$. By \cite[8.11]{cite:borel}, we may consider that $T$ splits over $\hat{k}$, some finite separable extension of $k$. 
    $T$ is contained in $C_{\mathcal{G}}(S)$, the Levi subgroup in our context 
    (Take $P^{\pm}=C_{\mathcal{G}}(S)U_{_{k}\Psi_{\pm}}$ as a pair of minimal opposite parabolic $k$-subgroups to recognize $C_{\mathcal{G}}(S)$ as the Levi subgroup, see \cite[20.4]{cite:borel}, \cite[20.6]{cite:borel}, and \cite[21.11]{cite:borel} for detail). 
    When we talk of a \textbf{homogeneous polynomial map of degree i} denoted by $f:V_{1}\to V_{2}$ for $R$-module's $V_{1}$ and $V_{2}$ (for arbitrary commutative ring $R$), we mean that the map $f$ has the property that $f(rv)=r^{i}f(v)$ for $v\in V_{1}$ and $r\in R$. 
    Furthermore, for any arbitrary $R$-algebra $A$, there is an induced map $A\otimes_{R} V_{1} \to A\otimes_{R} V_{2}$ by $a\otimes v \mapsto a^{i}\otimes f(v)$ that we denote by the same $f$ in an abuse of notation, this is consistent with the context and notations of \cite[sec: 6]{cite:nor_struc_of_iso_red_gp}. 
    We write 
    $\Psi:=\Phi(\mathcal{G},T)$ to be the absolute root system of $\mathcal{G}$, and 
    $\;_{k}\Psi:=\Phi(\mathcal{G},S)$ to be the relative root system of $\mathcal{G}$.

\end{note}

Petrov and Stavrova provided following generalization of absolute pinning isomorphism in the relative case (see \cite[thm: 2]{cite:ele_subgp_of_iso_red_gp}), we record the statements in our context: \bigskip 

\begin{defn}[Relative pinning maps]\label{n:Generalize unipotent element}
    For $a'\in\;_{k}\Psi$, take free $k$-module $V_{a'}$ as constructed in \cite[thm: 2]{cite:ele_subgp_of_iso_red_gp} (or \cite[lem: 6.1]{cite:nor_struc_of_iso_red_gp}). We have that $\hat{k}\otimes_{k}V_{a'}$ is free over $\hat{k}$. 
    We fix an ordered basis $M_{a'}:=\{e^{a'}_{\delta}\}_{\delta\in \eta(a')}$ of $\hat{k}\otimes_{k}V_{a'}$ over $\hat{k}$, and let $m_{a'}:=|M_{a'}|=|\eta(a')|$ be $\hat{k}\otimes_{k}V_{a'}$'s dimension over $\hat{k}$. 
    $T$ is a $\hat{R}:=\hat{k}[t,t^{-1}]$-split torus ($T$ is $\hat{k}$ split, so must be $\hat{R}$ split). 

    For any $u=\sum_{\delta\in \eta(a')} c_{\delta} \otimes e^{a'}_{\delta}\in \hat{R}\otimes_{k} V_{a'}$ (we sometimes skip writing the ``\;$\otimes$'' in the simple tensor elements if the context is clear), we extend the construction of absolute pinning isomorphisms to the \textbf{relative pinning maps} as following  
    \[x_{a'}(u):=\Bigl(\prod_{\delta\in \eta(a')} x_{\delta}(c_{\delta})\Bigr)\Bigl(\prod_{\theta\in \eta(2a')}x_{\theta}\bigl(p^{2}_{\theta}(u)\bigr)\Bigr)\;,\] 
    where for $a\in \Psi$, $x_{a}(\cdot)$ is the admissible isomorphism (i.e. absolute pinning isomorphism) as constructed in the split case (see \ref{d:Split reductive group and admissible isomorphism}). If for a root $*$, $x_{*}(\cdot)$ does not exist, we treat the specific $x_{*}(\cdot)$ as 1. 
    $p^{2}_{\theta}: \hat{R}\otimes_{k} V_{a'}\to \hat{R}\otimes k \cong \hat{R}$ is a homogeneous of degree 2 map induced by $p_{\theta}^{2}:V_{a'}\to k$. 

    For each relative root $a'$, we fix a basis $M^k_{a'}:=\{e^{a'}_1,\cdots, e^{a'}_{_{k}m_{a'}}\}$ of $V_{a'}$ over $k$ with $_{k}m_{a'}\geq 1$. We consider $V_{a'}$ embedded into $\hat{k}\otimes_{k}V_{a'}$ by $v\mapsto 1\otimes v$. 
    Then, $\displaystyle e^{a'}_{i}=\sum_{\delta\in\eta(a')}f_{\delta}e^{a'}_{\gamma}$ for $f_{\delta}\in\hat{k}$, hence, we have that 
    \[\displaystyle x_{a'}(ce^{a'}_{i})=x_{a'}\Bigl(\sum_{\delta\in\eta(a')}c\cdot f_{\delta}e^{a'}_{\gamma}\Bigr)=\Bigl(\prod_{\delta\in\eta(a')}x_{\delta}(c\cdot f_{\delta})\Bigr)\Bigl(\prod_{\theta\in\eta(2a')}x_{\theta}\bigl(p^{2}_{\theta}(\sum_{\delta\in\eta(a')}c\cdot f_{\delta}e^{a'}_{\gamma})\bigr)\Bigr)\;.\] 
    Notice that, since $_{k}\Psi$ is a root system, above construction did not need to consider root that is more than 2 multiple. While in case of \cite[lem: 6.1]{cite:nor_struc_of_iso_red_gp} and \cite[thm: 2]{cite:ele_subgp_of_iso_red_gp}, the corresponding sets of roots are not necessarily root systems. 
\end{defn}\bigskip 

\begin{note}[Properties of realtive pinning maps]\label{n:Properties of realtive pinning maps}
    Relative pinning maps come with many properties (see \cite[lem: 6.2]{cite:nor_struc_of_iso_red_gp} and \cite[lem: 7]{cite:ele_subgp_of_iso_red_gp} for reference), we record some of them: 
    \begin{enumerate}
        \item (\textbf{Additive property}) (\cite[lem: 6.2 (i)]{cite:nor_struc_of_iso_red_gp}): There exists a degree 2 homogeneous polynomial map $q^{2}_{a'}: V_{a'}\oplus V_{a'}\to V_{2a'}$, so that for a $k$-algebra $R$, and any $v,w\in R\otimes_{k} V_{a'}$, 
        we have that  
        \[x_{a'}(v)x_{a'}(w)=x_{a'}(v+w)x_{2a'}\Bigl(q^{2}_{a'}(v,w)\Bigr)\;.\] 
        
        \item (\textbf{Conjugation (by elements of Levi subgroup) property}) (\cite[lem:6.2 (ii)]{cite:nor_struc_of_iso_red_gp}): Recall that the Levi subgroup in our context is $C_{\mathcal{G}}(S)$.  
        We have for any $g\in C_{\mathcal{G}}(S)(k)$ there exist degree i homogeneous polynomial maps $\phi^{i}_{g,a'}:V_{a'}\to V_{ia'}$ for $i=1,2$, so that for a $k$-algebra R, and any $v\in R\otimes_{k} V_{a'}$,
        we have that  
        \[gx_{a'}(v)g^{-1}=x_{a'}\Bigl(\phi^{1}_{g,a'}(v)\Bigr)x_{2a'}\Bigl(\phi^{2}_{g,a'}(v)\Bigr)\;.\] 

        \item (\textbf{(Generalized) Chevalley commutator formula}) (\cite[lem:6.2 (iii)]{cite:nor_struc_of_iso_red_gp}): For $a',b'\in \;_{k}\Psi$ such that $ma'\neq -kb'$ for all $m,k\geq 1$, there exist polynomial map(s)  
        \[N_{a',b';i,j}:V_{a'}\times V_{b'}\to V_{ia'+jb'}\] 
        for $i,j>0$. The above map $N_{a',b';i,j}$ is homogeneous of degree i in the first variable, and of degree j in the second variable (implying that $(r_{1}\otimes v, r_{2}\otimes w)\mapsto r_{1}^{i}r_{2}^{j}\otimes N_{a',b';i,j}(v,w)$ for the induced map). Furthermore, for a $k$-algebra $R$, any $v\in R\otimes_{k}V_{a'}$, and any $w\in R\otimes_{k} V_{b'}$, 
        we have that 
        \[\bigl[x_{a'}(v),x_{b'}(w)\bigr]=\prod_{i,j>0, ia'+jb'\in \;_{k}\Psi} x_{ia'+jb'}\bigl(N_{a',b';i,j}(v,w)\bigr)\;.\] 

        \item (\textbf{Generating property}) (\cite[lem: 7]{cite:ele_subgp_of_iso_red_gp}): Let $_{k}\Psi'$ be a \textbf{unipotent closed subset} of $_{k}\Psi$ (that is, a subset of $_{k}\Psi$ that is closed under addition when the sum is also a k-root, and does not contain collinear oppositely directed $k$-roots.). Take $U_{_{k}\Psi'}:=\langle U_{a'}|a'\in _{k}\Psi' \rangle$, then $U_{_{k}\Psi'}(R)$ is generated by the collection of $x_{a'}(r\otimes e^{a'}_{i})$'s for $r\in R$, $a\in \;_{k}\Psi'$, and $e^{a'}_{i}\in M^{k}_{a'}$.
        
        As a further comment: The notion of ``unipotent closed set'' is same as the notion of ``special set'' in the context of \cite[14.5]{cite:borel}, see \cite[14.7]{cite:borel} for detail.
    \end{enumerate}
\end{note}\bigskip

%% file: sections/Main_theorem_and_its_proof.tex
We maintain the context stated in \ref{n:Context in needed results}. We write $R:=k[t,t^{-1}]$. Furthermore, we require the relative root system be irreducible ($\mathcal{G}$ being almost simple guarantees this requirement). We abuse notation and write $\mathcal{G}(\overline{k(t)})$ to be the same as $\mathcal{G}$. 

As a reminder: As in \cite[15.3.6]{cite:TAS_LAG}, we understand the group of cocharacters of $S$ to be a subgroup of group of cocharacters of $T$, and we take the unique coroot for relative root in the similar fashion akin to the absolute root case (see \cite[7.1.8]{cite:TAS_LAG}). More details about the relative root system $_{k}\Psi$ is stated in \cite[sec: 15.3]{cite:TAS_LAG}. 

In this section, we first construct the affine root groups, and make observation on their properties in subsection \ref{subsec:Preparations}. Then, we provide the proof of the main theorem in subsection \ref{subsec:Main theorem}.

\subsection{Preparations}\label{subsec:Preparations} 

\begin{defn}[Construction of affine root groups]\label{n:Construction of Affine root group in non-split case}
Take $a' \in \;_{k}\Psi$, and $l\in\mathbb{Z}$. We define the \textbf{affine root groups in the not necessarily split case} as
\[U_{\alpha_{a',l}}:=\Bigl\langle x_{a'}(ct^{-l}\otimes e^{a'}_{i}),x_{2a'}(ct^{-2l}\otimes e^{2a'}_{j}) \Big| c\in k, e^{a'}_{i}\in M^{k}_{a'}, e^{2a'}_{j}\in M^{k}_{2a'} \Bigr\rangle\;. \] 
In context of the not necessarily split isotropic reductive groups, we denote $\Phi$ to be the set of affine roots associated with $_{k}\Psi$. 
We also make use of the extended case where $2l$ is within full and half integers, where the construction stays the same (although, technically, the constructed groups are no longer affine root groups). 

Note that, in the case of $(a')=\{a',2a'\}$, we no longer have $U_{\alpha_{2a',l}}\leq U_{\alpha_{a',l}}$. But we still have $U_{\alpha_{a',l}}\leq U_{a'}$ for $a'\in \;_{k}\Psi$. An alternative but equivalent construction of the affine root groups is presented at \ref{n:Alternative form of affine root group in non-splt case}. 
\end{defn}\bigskip

\begin{obs}[Commutativity within the same affine root group]\label{o:Commutativity within the same affine root group}
    Take $r\in R$, and $\delta\in \eta(a')$, then we have $x_{\delta}(r)\in U_{a'}$. Take $s\in R$, and $\gamma\in \eta(2a')$, then we have $x_{\gamma}(s)\in U_{2a'}=Z(U_{a'})$. This means that 
            \[x_{a'}(v)\; \text{and} \;x_{2a'}(w)\; \text{commute}\] 
            for all $v\in R\otimes_{k} V_{a'}$, and $w\in R\otimes_{k} V_{2a'}$. 
\end{obs}\bigskip

\begin{obs}[Consequences of properties of relative pinning maps]\label{n:Important proteries of generalize unipotent elements and affine root group}
    We apply properties of relative pinning maps: 
    \begin{enumerate}
        \item[(a)] Consequence of the additive property:  

        \begin{itemize}
            \item Take any $k$-root $a'$ such that $(a')=\{a'\}$. Take $c,d\in k$, and (not necessarily distinct elements) $e^{a'}_{i},e^{a'}_{j}\in M^k_{a'}$,  
            we have 
            \[x_{a'}\Bigl(\sum_{i}c_{i}t^{-l}\otimes e^{a'}_{i}\Bigr)=\prod_{i}x_{a'}\bigl(c_{i}t^{-l}\otimes e^{a'}_{i}\bigr)\in U_{\alpha_{a',l}}\;.\] 
            Together with the fact that $x_{a'}(0)=1$ (see \cite[thm: 2]{cite:ele_subgp_of_iso_red_gp}), above implies that $x_{a'}(v)^{-1}=x_{a'}(-v)$ for all $v\in R\otimes_{k} V_{a'}$. 
            
            \item Take any $k$-root $a'$ such that $(a')=\{a',2a'\}$. Take $c,d\in k$, and (not necessarily distinct elements) $e^{a'}_{i},e^{a'}_{j}\in M^{k}_{a'}$, we have  
            \[x_{a'}\bigl(ct^{-l}\otimes e^{a'}_{i}+dt^{-l}\otimes e^{a'}_{j}\bigr)
            =x_{a'}\bigl(ct^{-l}\otimes e^{a'}_{i}\bigr)x_{a'}\bigl(dt^{-l}\otimes e^{a'}_{j}\bigr) x_{2a'}\Bigl(-q^{2}_{a'}\bigl(ct^{-l}\otimes e^{a'}_{i},dt^{-l}\otimes e^{a'}_{j}\bigr)\Bigr)\;.\] 
            Notice, also,
            \[q^{2}_{a'}\bigl(ct^{-l}\otimes e^{a'}_{i},dt^{-l}\otimes e^{a'}_{j}\bigr)=(t^{-2l})\otimes q^{2}_{a'}\bigl(ce^{a'}_{i},de^{a'}_{j}\bigr)
            =\sum_{h=1}^{_{k}m_{2a'}} (f_{h})(t^{-2l}) \otimes e^{2a'}_{h}\] 
            for $f_{h}\in k$. To verify that $f_{h}$ is contained in $k$, recall that $e^{a'}_{i},e^{a'}_{j},e^{a'}_{h}$ are contained in $M^{k}_{a'}$, where $M^{k}_{a'}$ is considered a generating set of $V_{a'}$ over $k$. 
            Together with the result in first bullet point, we infer that 
            {\footnotesize\[x_{a'}\bigl(ct^{-l}\otimes e^{a'}_{i}+dt^{-l}\otimes e^{a'}_{j}\bigr)=x_{a'}\bigl(ct^{-l}\otimes e^{a'}_{i}\bigr)x_{a'}\bigl(dt^{-l}\otimes e^{a'}_{j}\bigr) x_{2a'}\Bigl(\sum_{h=1}^{_{k}m_{2a'}} -(f_{h})(t^{-2l}) \otimes e^{2a'}_{h}\Bigr)\in U_{\alpha_{a',l}}\label{eq:prodtosum} \tag{*}\;,\]}
            which implies, through induction of changing one of the simple tensors into a sum of finite many simple tensors, 
            \[x_{a'}\Bigl(\sum_{i}c_{i}t^{-l}\otimes e^{a'}_{i}\Bigr)\in U_{\alpha_{a',l}}\;.\]

            \item In the case where $a'\in\; _{k}\Psi$ so that $2a'\notin\; _{k}\Psi$, we have 
            \[x_{a'}(v)^{-1}=x_{a'}(-v)\;.\] 
            In the case where $a',2a'\in\; _{k}\Psi$, we have 
            \[x_{a'}(v)x_{a'}(-v)=x_{2a'}\Bigl(q^{2}_{a'}(v,-v)\Bigr)\;,\] 
            and, so, 
            \[x_{a'}(v)^{-1}=x_{a'}(-v)x_{2a'}\Bigl(-q^{2}_{a'}(v,-v)\Bigr)\;.\] 

        \end{itemize}

        \item[(b)] Consequence of the conjugation property: We show that \textbf{$C_{\mathcal{G}}(S)(k)$ normalizes each affine root group in the not necessarily split case}, we only need to consider the generators. Take $g\in C_{\mathcal{G}}(S)(k)$ and consider following cases: 

        \begin{itemize}
            \item If $(a')=\{a'\}$, we have 
            \[gx_{a'}\bigl(ct^{-l}\otimes e^{a'}_{i}\bigr)g^{-1}=x_{a'}\Bigl(\phi^{1}_{g,a'}(ct^{-l}\otimes e^{a'}_{i})\Bigr)=x_{a'}\Bigl((ct^{-l})\otimes \phi^{1}_{g,a'}(e^{a'}_{i})\Bigr)\;.\] 
            We take $\phi^{1}_{g,a'}(e^{a'}_{i})=\sum_{h=1}^{_{k}m_{a'}} f_{h} e^{a'}_{h}$ for $f_{h}\in k$, and obtain  
            \[gx_{a'}\bigl(ct^{-l}\otimes e^{a'}_{i}\bigr)g^{-1}=x_{a'}\Bigl((ct^{-l})\otimes \sum_{h=1}^{_{k}m_{a'}} f_{h} e^{a'}_{h}\Bigr)=x_{a'}\Bigl(\sum_{h=1}^{_{k}m_{a'}} (cf_{h})(t^{-l})\otimes e^{a'}_{h}\Bigr)\in U_{\alpha_{a',l}}\;.\] 
            
            \item If $(a')=\{a',2a'\}$, we have  
            \[gx_{a'}\bigl(ct^{-l}\otimes e^{a'}_{i}\bigr)g^{-1}=x_{a'}\Bigl(\phi^{1}_{g,a'}(ct^{-l}\otimes e^{a'}_{i})\Bigr)x_{2a'}\Bigl(\phi^{2}_{g,a'}(ct^{-l}\otimes e^{a'}_{i})\Bigr)\;.\]
            We take $\phi^{1}_{g,a'}(e^{a'}_{i}) =\sum_{h=1}^{_{k}m_{a'}} f^{1}_{h} e^{a'}_{h}$ and $\phi^{2}_{g,a'}(e^{a'}_{i})=\sum_{h=1}^{_{k}m_{2a'}} f^{2}_{h} e^{2a'}_{h}$ for $f^{1}_{h},f^{2}_{h}\in k$, and obtain 
            \[gx_{a'}\bigl(ct^{-l}\otimes e^{a'}_{i}\bigr)g^{-1}=x_{a'}\Bigl((ct^{-l})\otimes \sum_{h=1}^{_{k}m_{a'}} f^{1}_{h} e^{a'}_{h}\Bigr)x_{2a'}\Bigl((c^{2}t^{-2l})\otimes \sum_{h=1}^{_{k}m_{2a'}} f^{2}_{h} e^{2a'}_{h}\Bigr)\] 
            \[=x_{a'}\Bigl(\sum_{h=1}^{_{k}m_{a'}} (cf^{1}_{h})(t^{-l})\otimes e^{a'}_{h}\Bigr)x_{2a'}\Bigl(\sum_{h=1}^{_{k}m_{2a'}} (c^{2}f^{2}_{h})(t^{-2l})\otimes e^{2a'}_{h}\Bigr)\in U_{\alpha_{a',l}}\;.\] 
            Notice that $gx_{2a'}(ct^{-2l}\otimes e^{2a'}_{i})g^{-1}\in U_{\alpha_{a',l}}$ is already proven in the previous bullet point. 
        \end{itemize}

        \item[(c)] The Chevalley commutator formula is a key tool for verifying (RGD1), we detail this later within the proof of the main theorem.

        \item[(d)] Consequence of the generating property: Take $a'\in\;_{k}\Psi$, then $(a')$ is an unipotent closed subset, and, so, we have that $U_{(a')}(R)$ is generated by the collection of $x_{a'}(r\otimes e^{a'}_{i})$'s for $r\in R$, $a\in (a')$, and $e^{a'}_{i}\in M^k_{a'}$. Hence, we have $U_{\alpha_{a',0}}=U_{a'}(k)$ as a consequence of the generating property. 
        
        Using $_{k}\Psi_{\pm}$ as an unipotent closed subset of $_{k}\Psi$, we infer that $U_{_{k}\Psi_{\pm}}(R)$ is generated by the collection of $x_{a'}(r\otimes e^{a'}_{i})$'s for $r\in R$, $a\in \;_{k}\Psi_{\pm}$, and $e^{a'}_{i}\in M^k_{a'}$. 
        As a consequence, the \textbf{elementary subgroup} $\mathcal{G}(R)^{+}:=\langle U_{_{k}\Psi_{+}}(R),U_{_{k}\Psi_{-}}(R) \rangle $ is generated by the collection of $x_{a'}(r\otimes e^{a'}_{i})$'s for $r\in R$, $a\in \;_{k}\Psi$, and $e^{a'}_{i}\in M^k_{a'}$.

    \end{enumerate}
\end{obs}\bigskip

\begin{obs}[Alternative form of affine root groups in the not necessarily splt case]\label{n:Alternative form of affine root group in non-splt case}
    \[U_{\alpha_{a',l}}:=\Bigl\{x_{a'}\bigl(\sum_{i=1}^{_{k}m_{a'}}c_{i}t^{-l}\otimes e^{a'}_{i}\bigr)x_{2a'}\bigl(\sum_{i=1}^{_{k}m_{2a'}}d_{i}t^{-2l}\otimes e^{2a'}_{i}\bigr) \Big| c_{i},d_{i}\in k\Bigr\}\] 
    is an alternative and equivalent way to construct the affine root groups. 
    \begin{Sketch}
        This observation is largely a consequence of \ref{n:Important proteries of generalize unipotent elements and affine root group} (a): We start with an arbitrary element in $U_{\alpha_{a',l}}$ represented as products of generators in \ref{n:Construction of Affine root group in non-split case}, for instance, take $h_{i},w_{j}\in k$ and consider the element 
        \[\Bigl(\prod_{i} x_{a'}(h_{i} t^{-l}\otimes e^{a'}_{i})\Bigr)\Bigl(\prod_{j} x_{2a'}(w_{i} t^{-2l}\otimes e^{2a'}_{j})\Bigr)\;.\] 
        Notice, in above, we already moved all the $x_{a'}(\cdot)$ terms to the left, and all the $x_{2a'}(\cdot)$ terms to the right (we can do this because of \ref{o:Commutativity within the same affine root group}). 
        Consulting line \eqref{eq:prodtosum} from (a) of \ref{n:Important proteries of generalize unipotent elements and affine root group}, 
         we obtain 
        \[\prod_{i} x_{a'}(h_{i} t^{-l}\otimes e^{a'}_{i})=x_{a'}\bigl(\sum_{i} h_{i} t^{-l}\otimes e^{a'}_{i}\bigr)x_{2a'}\bigl(\sum_{j=1}^{_{k}m_{2a'}} f_{j}t^{-2l}\otimes e^{2a'}_{j}\bigr)\] 
        for $f_{j}$'s in $k$. 
        Above is the general form of elements generated only by the $x_{a'}(\cdot)$ terms. 
        Then, consulting first bullet point of \ref{n:Important proteries of generalize unipotent elements and affine root group} (a), 
        we multiply in all the $x_{2a'}(\cdot)$ terms, and 
        obtain that any element of $U_{\alpha_{a',l}}$ has the form 
        \[x_{a'}\bigl(\sum_{i=1}^{_{k}m_{a'}}c_{i}t^{-l}\otimes e^{a'}_{i}\bigr)x_{2a'}\bigl(\sum_{i=1}^{_{k}m_{2a'}}d_{i}t^{-2l}\otimes e^{2a'}_{i}\bigr)\]
        for $c_{i},d_{i}\in k$. In other words, we have shown that 
        \[U_{\alpha_{a',l}}\subset\Bigl\{x_{a'}\bigl(\sum_{i=1}^{_{k}m_{a'}}c_{i}t^{-l}\otimes e^{a'}_{i}\bigr)x_{2a'}\bigl(\sum_{i=1}^{_{k}m_{2a'}}d_{i}t^{-2l}\otimes e^{2a'}_{i}\bigr) \Big| c_{i},d_{i}\in k\Bigr\}\;.\]
        To, further, replace $\subset$ with $=$, we utilize the first two bullet points in (a) of \ref{n:Important proteries of generalize unipotent elements and affine root group}. 
    \end{Sketch}
\end{obs}\bigskip

\begin{lem}[Relative pinning maps and affine root groups conjugated by elements in the image of coroots]\label{o:Relative pinning maps conjugated by cocharacters}
    
    We have, for $n,l$ in full and half integers,  
    \[a'^{\vee}(t^{-l/2})U_{\alpha_{b',n}}a'^{\vee}(t^{-l/2})^{-1}= U_{\alpha_{b',\frac{l<b',a'>}{2}+n}} \label{eq:coroot-move-afrg}\tag{!}\;.\]

    \begin{Proof}
        We focus on the generators of the group 
        \[U_{\alpha_{b',n}}:=\Bigl\langle x_{b'}(ct^{-n}\otimes e^{b'}_{i}),x_{2b'}(ct^{-2n}\otimes e^{2b'}_{j}) \Big| c\in k, e^{b'}_{i}\in M^k_{b'}, e^{2b'}_{j}\in M^k_{2b'} \Bigr\rangle\;.\] 
        When we talk about affine root groups, we consider $n\in \mathbb{Z}$; but in the following, we consider n to be half and full integers (so that $U_{\alpha_{b',n}}\subset \mathcal{G}(k[t^{-1/2},t^{1/2}])$). One may find item 1 of \ref{n:Set up for non-split case} helpful for the following. 
        Consider $l$ in full and half integers (the following calculations are made within $\mathcal{G}(k[t^{1/4},t^{-1/4}])$) 

        \begin{itemize}
            \item In the case of relative root $2b'$, we start with the element 
            \[a'^{\vee}(t^{-l/2})x_{2b'}(ct^{-2n}\otimes e^{2b'}_{i})a'^{\vee}(t^{-l/2})^{-1}\;.\] 
            We write $\displaystyle e^{2b'}_{i}=\sum_{\theta\in \eta(2b')}h_{\theta}e_{\theta}^{a'}$ for $h_{\theta}\in \hat{k}$, and, hence, have that $\displaystyle x_{2b'}(ce^{2b'}_{i}) = \prod_{\theta\in \eta(2b')} x_{\theta}(c\cdot h_{\theta})$. The element we started with then becomes 
            \[a'^{\vee}(t^{-l/2})x_{2b'}(ct^{-2n}\otimes e^{2b'}_{i})a'^{\vee}(t^{-l/2})^{-1} 
            =a'^{\vee}(t^{-l/2})\Bigl(\prod_{\theta\in \eta(2b')} x_{\theta}(ct^{-2n}\cdot h_{\theta})\Bigr)a'^{\vee}(t^{-l/2})^{-1}\] 
            \[= \prod_{\theta\in \eta(2a')} x_{\theta}\bigl(ct^{-l<b',a'>-2n}\cdot h_{\theta}\bigr) = x_{2b'}\bigl(ct^{-l<b',a'>-2n}\otimes e^{2b'}_{i}\bigr)\in U_{\alpha_{b',n+\frac{l<b',a'>}{2}}}\;.\]

            \item In the case of relative root $b'$, we start with the element  
            \[a'^{\vee}(t^{-l/2})x_{b'}(ct^{-n}\otimes e^{b'}_{i})a'^{\vee}(t^{-l/2})^{-1}\;.\]
            We write $\displaystyle e^{b'}_{i}=\sum_{\delta\in\eta(b')}f_{\delta}e^{b'}_{\gamma}$ for $f_{\delta}\in\hat{k}$, then, we have 
            \[x_{b'}(ce^{b'}_{i})=\Bigl(\prod_{\delta\in\eta(b')}x_{\delta}(c\cdot f_{\delta})\Bigr)\Bigl(\prod_{\theta\in\eta(2b')}x_{\theta}\bigl(p^{2}_{\theta}(\sum_{\delta\in\eta(b')}c\cdot f_{\delta}e^{b'}_{\gamma})\bigr)\Bigr)\;.\] We infer 
            \[a'^{\vee}(t^{-l/2})x_{b'}(ct^{-n}\otimes e^{b'}_{i})a'^{\vee}(t^{-l/2})^{-1}\]
            \[=a'^{\vee}(t^{-l/2})\Bigl(\prod_{\delta\in\eta(b')}x_{\delta}(ct^{-n}\cdot f_{\delta})\Bigr)\Bigl(\prod_{\theta\in\eta(2b')}x_{\theta}\bigl(p^{2}_{\theta}(\sum_{\delta\in\eta(b')}ct^{-n}\cdot f_{\delta}e^{b'}_{\gamma})\bigr)\Bigr)a'^{\vee}(t^{-l/2})^{-1}\;.\] 
            Notice that we can insert $a'^{\vee}(t^{-l/2})a'^{\vee}(t^{-l/2})^{-1}$ in between each term of the absolute pinning isomorphisms, with the conjugation equality for a single absolute pinning isomorphism associated with the absolute root $\delta$ that projects to $b'$ as 
            \[a'^{\vee}(t^{-l/2})x_{\delta}(ct^{-n}\cdot f_{\delta})a'^{\vee}(t^{-l/2})^{-1}=x_{\delta}\bigl(ct^{-\frac{l<b',a'>}{2}-n}\cdot f_{\delta}\bigr)\;.\] 
            For the $k$-module $V_{b'}$, we have that $p^{2}_{\theta}: \hat{k}[t^{1/2},t^{-1/2}]\otimes_{k}V_{b'}\to \hat{k}[t^{1/2},t^{-1/2}]\otimes_{k} k\cong \hat{k}[t^{1/2},t^{-1/2}]$ is induced by homogeneous degree 2 map $p^{2}_{\theta}: V_{b'}\to k$ by $r\otimes v\mapsto r^{2}\otimes p^{2}_{\theta}(v)$. We may write $\displaystyle p^{2}_{\theta}(\sum_{\delta\in\eta(b')}c\cdot f_{\delta}e^{b'}_{\gamma})=k_{\theta}$, and, hence, $\displaystyle p^{2}_{\theta}(\sum_{\delta\in\eta(b')}ct^{-m}\cdot f_{\delta}e^{b'}_{\gamma})=p^{2}_{\theta}(t^{-m}\sum_{\delta\in\eta(b')}c\cdot f_{\delta}e^{b'}_{\gamma}) =k_{\theta}t^{-2m}$ for $2m$ in half and full integers.  
            We obtain the conjugation equality for a single absolute pinning isomorphism associated with the absolute root $\theta$ that projects to $2b'$ 
            {\footnotesize\[a'^{\vee}(t^{-l/2})x_{\theta}\Bigl(p^{2}_{\theta}\bigl(\sum_{\delta\in\eta(b')}ct^{-n}\cdot f_{\delta}e^{b'}_{\gamma}\bigr)\Bigr)a'^{\vee}(t^{-l/2})^{-1}=a'^{\vee}(t^{-l/2})x_{\theta}(k_{\theta}t^{-2n})a'^{\vee}(t^{-l/2})^{-1}
            =x_{\theta}\bigl(k_{\theta}t^{-l<b',a'>-2n}\bigr)\;.\] }
            Then, the element we started with becomes 
            {\small\[a'^{\vee}(t^{-l/2})x_{b'}(ct^{-n}\otimes e^{b'}_{i})a'^{\vee}(t^{-l/2})^{-1}
            =\Bigl(\prod_{\delta\in\eta(b')} x_{\delta}\bigl(ct^{-\frac{l<b',a'>}{2}-n}\cdot f_{\delta}\bigr)\Bigr)\Bigl(\prod_{\theta\in\eta(2b')}x_{\theta}\bigl(k_{\theta}t^{-l<b',a'>-2n}\bigr)\Bigr)\] 
            \[=\Bigl(\prod_{\delta\in\eta(b')} x_{\delta}\bigl(ct^{-\frac{l<b',a'>}{2}-n}\cdot f_{\delta}\bigr)\Bigr)\Bigl(\prod_{\theta\in\eta(2b')}x_{\theta}\bigl(p^{2}_{\theta}(\sum_{\delta\in\eta(b')} ct^{-\frac{l<b',a'>}{2}-n}f_{\delta}e^{b'}_{\delta})\bigr)\Bigr)\] 
            \[=x_{b'}(ct^{-\frac{l<b',a'>}{2}-n}\otimes e^{b'}_{i})\in U_{\alpha_{b',\frac{l<b',a'>}{2}+n}}\;.\] }
        \end{itemize}
        Above shows $a'^{\vee}(t^{-l/2})U_{\alpha_{b',n}}a'^{\vee}(t^{-l/2})^{-1}\subset U_{\alpha_{b',\frac{l<b',a'>}{2}+n}}$. On the other hand, notice that $a'^{\vee}(t^{l/2})U_{\alpha_{b',\frac{l<b',a'>}{2}+n}}a'^{\vee}(t^{-l/2})\subset U_{\alpha_{b',\frac{l<b',a'>}{2}+n-\frac{l<b',a'>}{2}}}=U_{\alpha_{b',n}}$. 
    \end{Proof}
\end{lem}\bigskip

\begin{cor}[Isomorphism between affine root groups by coroots]\label{o:Isomorphism between affine root groups by coroots}
    By equation \eqref{eq:coroot-move-afrg} and the fact that conjugation by a group element is a group isomorphism onto its image, we see that any element of $x\in U_{\alpha_{a'},l}$ (resp. more generally, $x\in U_{\alpha_{b',\frac{l<b',a'>}{2}+n}}$) can be expressed as 
    \[x=a'^{\vee}(t^{-l/2})ua'^{\vee}(t^{-l/2})^{-1}\] 
    for some unique $u\in U_{a'}(k)$ (resp. more generally, $u\in U_{\alpha_{b',n}}$). 
\end{cor}\bigskip

\begin{lem}[Preparation for (RGD2) in the not necessarily split case]\label{n:Preparation for (RGD2) in non-split case} 
    Take an affine root $\alpha:=\alpha_{b',m}$. By \ref{o:Isomorphism between affine root groups by coroots}, for any $x\in U_{\alpha}$, there is an unique $u\in U_{b'}(k)$ so that $\displaystyle x=b'^{\vee}(t^{-l/2})ub'^{\vee}(t^{-l/2})^{-1}$, we take and fix said elements $x$ and $u$. There exists an element $w_{a',u}$ so that the element $\displaystyle  w_{a',x}(t^{-l}):=a'^{\vee}(t^{-l/2})w_{a',u}a'^{\vee}(t^{-l/2})^{-1}$ is contained in $\displaystyle U_{-\alpha_{a',l}}xU_{-\alpha_{a',l}}$, and 
    \[w_{a',x}(t^{-l})U_{\alpha_{b',m}}w_{a',x}(t^{-l})^{-1}=U_{s_{a',l}(\alpha_{b',m})}\;.\] 
    
    \begin{Proof}

        We, first, provide an informal summary: Within the action of the element $w_{a',x}(t^{-l})$ on the affine root group $U_{\alpha}$, the element $w_{a',u}$ ``acts'' only on the index $b'$ by reflection in respect to $a'$ (resulting in $s_{a'}(b')$); and, $l$, being the power of $t$, determines the change in the integer coordinate of $\alpha$ while not affecting the coordinate of the relative root. 

        \begin{enumerate}
            \item[\textbf{(Step 1)}] (Finding the element $w_{a',u}$): We take, for $a'\in \;_{k}\Psi$, an element $n_{a'}=n_{-a'}\in N_{\mathcal{G}}(S)(k)$ as in c) of \cite[Satz: 25]{cite:Reducktive_Gruppen_AP} (the element originates in \cite[21.2]{cite:borel}, and is also recorded in \cite[15.3.5]{cite:TAS_LAG}). This element $n_{-a'}$ enjoys following properties (DR4), and (DR5) (see \cite[Satz: 27]{cite:Reducktive_Gruppen_AP}):  
            \begin{enumerate}
                \item[(DR4)] For any $a'\in \; _{k}\Psi$: \[U_{a'}(k)\setminus\{e\}\subset U_{-a'}(k)C_{\mathcal{G}}(S)(k)n_{-a'}U_{-a'}(k)\;.\] 
                \item[(DR5)] For any $a',b'\in\; _{k}\Psi$, $w\in C_{\mathcal{G}}(S)(k)n_{-a'}$: \[wU_{b}(k)w^{-1}=U_{s_{a'}(b')}(k)\;.\] 
            \end{enumerate} 
            We turn our sight to \cite[21.1]{cite:borel}, where it is stated that $N_{\mathcal{G}}(S)/C_{\mathcal{G}}(S)$ acts on $S$ and induces an action on the cocharacter group $X_{*}(S)$. Take $n\in N_{\mathcal{G}}(S)$, $\chi\in X_{*}(S)$, and $\lambda\in \overline{k(t)}$,  
            \[nC_{\mathcal{G}}(S)\cdot \chi (\lambda):=n\chi(\lambda)n^{-1}\] 
            should describe this action. 
            This action extends to an action on $X_{*}(S)\otimes_{\mathbb{Z}}\mathbb{R}$. 
            $n_{-a'}=n_{a'} \in N_{\mathcal{G}}(S)(k)$ is taken so that the element $n_{-a'}C_{\mathcal{G}}(S)\in N_{\mathcal{G}}(S)/C_{\mathcal{G}}(S)$ induces an orthogonal transformation on $X_{*}(S)\otimes_{\mathbb{Z}}\mathbb{R}$ while fixing the hyperplane $X_{*}(S_{a'})\otimes_{\mathbb{Z}}\mathbb{R}$ (where $S_{a'}$ is the identity component of $Ker(a')$, $n_{-a'}$'s counterpart is $r_{a'}$ in \cite[21.2]{cite:borel}, or $s_{a'^{\vee}}$ in \cite[15.3.5-15.3.8]{cite:TAS_LAG}). 
            See about reflections in $X_{*}(S)\otimes_{\mathbb{Z}}\mathbb{R}$ at \ref{n:Reduced root system in reductive group}, where the formula for reflections in the co-space ``stay the same'' when generalized to the not necessarily split case. In particular, $s_{a'}(b')^{\vee}=s_{a'^{\vee}}(b'^{\vee})$ is still a consequence of the formula for reflections in the co-space and the bilinearity of $<\cdot,\cdot>$ when understood as the composition of characters and cocharacters. 
            To sum up, this means that $n_{-a'}$ is taken as the reflection in respect to the hyperplane $X_{*}(S_{a'})\otimes_{\mathbb{Z}}\mathbb{R}$, hence, satisfies the condition 
            \[n_{-a'}b'^{\vee}(\lambda)n_{-a'}^{-1}=s_{a'}(b')^{\vee}(\lambda) \tag{***}\label{eq:on na}\]
            for all $\lambda\in \overline{k(t)}$. By (DR4), take any $u\in U_{a'}(k)^{*}$, we have that 
            \[u\in U_{-a'}(k)C_{\mathcal{G}}(S)(k)n_{-a'}U_{-a'}(k)\;.\] 
            We take and fix $w_{a',u}\in C_{\mathcal{G}}(S)(k)n_{-a'}$, so that we have $w_{a',u}\in U_{-a'}(k)uU_{-a'}(k)$. We denote $w_{a',u}=h_{a',u}n_{-a'}$ for $h_{a',u}\in C_{\mathcal{G}}(S)(k)$. And by (DR5), this choice of $w_{a',u}$ enjoys the property 
            \[w_{a',u}U_{b'}(k)w_{a',u}^{-1}=U_{s_{a'}(b')}(k)\;.\] 
            Now we extend the property brought to us by (DR5) with condition \eqref{eq:on na}: 
            \[w_{a',u}U_{\alpha_{b',n}}w_{a',u}^{-1}=h_{a',u}n_{-a'}U_{\alpha_{b',n}}n_{-a'}^{-1}h_{a',u}^{-1}=h_{a',u}n_{-a'}b'^{\vee}(t^{-n/2})U_{\alpha_{b',0}}b'^{\vee}(t^{-n/2})^{-1}n_{-a'}^{-1}h_{a',u}^{-1}\] 
            \[=h_{a',u}s_{a'}(b')^{\vee}(\lambda t^{-n/2})n_{-a'}U_{b'}(k)n_{-a'}^{-1}s_{a'}(b')^{\vee}(\lambda t^{-n/2})^{-1}h_{a',u}^{-1}\;.\] 
            Keep in mind that $U_{\alpha_{b',0}}=U_{b'}(k)$
            when verifying above statement. We continue the calculation:  
            \[w_{a',u}U_{\alpha_{b',n}}w_{a',u}^{-1}=h_{a',u}s_{a'}(b')^{\vee}(\lambda t^{-n/2})U_{s_{a'}(b')}(k)s_{a'}(b')^{\vee}(\lambda t^{-n/2})^{-1}h_{a',u}^{-1}\] 
            \[=h_{a',u}U_{\alpha_{s_{a'}(b'),n}}h_{a',u}^{-1}\subset U_{\alpha_{s_{a'}(b'),n}}\;,\] 
            where the last inclusion is by the conjugation property.
            We have shown that 
            \[w_{a',u}U_{\alpha_{b',n}}w_{a',u}^{-1}\subset U_{\alpha_{s_{a'}(b'),n}}\;.\] 
            With similar steps, we obtain 
            $w_{a',u}^{-1}U_{\alpha_{s_{a'}(b'),n}}w_{a',u}\subset U_{\alpha_{b',n}}$,  
            which infers \[w_{a',u}U_{\alpha_{b',n}}w_{a',u}^{-1}= U_{\alpha_{s_{a'}(b'),n}}\;.\] 
            Notice also, it is by construction that  
            \[w_{a',u_{1}}w_{a',u_{2}}^{-1}\in C_{\mathcal{G}}(S)(k)\;.\] 
            To sum up, we now have found, for all $a',b'\in \;_{k}\Psi$, and non-identity element $u\in U_{a'}(k)=U_{\alpha_{a',0}}$, an element $w_{a',u}\in U_{-a'}(k)uU_{-a'}(k)\subset \mathcal{G}(k)\subset \mathcal{G}(k[t,t^{-1}])$, for $\forall n\in\mathbb{Z}/2$, 
            so that 
            \[w_{a',u}U_{\alpha_{b',n}}w_{a',u}^{-1}=U_{\alpha_{s_{a'}(b'),n}}\subset \mathcal{G}(k[t^{-1/2},t^{1/2}])\text{, and } w_{a',u_{1}}w_{a',u_{2}}^{-1}\in C_{\mathcal{G}}(S)(k)\] 
            for all $u_{1},u_{2}\in U_{a'}(k)^{*}$. 
            \item[\textbf{(Step 2)}] (Constructing the element $w_{a',x}(t^{-l})$): By using the element $w_{a',u}$ in above and \ref{o:Isomorphism between affine root groups by coroots}, we can define and see that for any non identity element $x=a'^{\vee}(t^{-l/2})ua'^{\vee}(t^{-l/2})^{-1}\in U_{\alpha_{a',l}}$ where $u\in U_{a'}(k)$:  
            \[w_{a',x}(t^{-l}):=a'^{\vee}(t^{-l/2})w_{a',u}a'^{\vee}(t^{-l/2})^{-1}=a'^{\vee}(t^{-l/2})w_{a',u}a'^{\vee}(t^{l/2})\]
            \[\in a'^{\vee}(t^{-l/2})U_{-a'}(k)uU_{-a'}(k)a'^{\vee}(t^{l/2})\] 
            \[ = a'^{\vee}(t^{-l/2})U_{-a'}(k)a'^{\vee}(t^{l/2})\Bigl(a'^{\vee}(t^{-l/2})ua'^{\vee}(t^{-l/2})^{-1}\Bigr)a'^{\vee}(t^{-l/2})U_{-a'}(k)a'^{\vee}(t^{l/2})\] 
            \[ = U_{-\alpha_{a',l}}\Bigl(a'^{\vee}(t^{-l/2})ua'^{\vee}(t^{-l/2})^{-1}\Bigr)U_{-\alpha_{a',l}}=U_{-\alpha_{a',l}}xU_{-\alpha_{a',l}}\;.\] 
            We claim that the element $w_{a',x}(t^{-l})$ has the property that for $m\in\mathbb{Z}$, 
            \[w_{a',x}(t^{-l})U_{\alpha_{b',m}}w_{a',x}(t^{-l})^{-1}=U_{\alpha_{s_{a'}(b'),m-l<b',a'>}}=U_{s_{a',l}(\alpha_{b',m})}\;.\] 
            To verify this claim, we first check 
            \[a'^{\vee}(t^{l/2})U_{\alpha_{b',m}}a'^{\vee}(t^{l/2})^{-1}= U_{\alpha_{b',\frac{-l<b',a'>}{2}+m}}\;,\] 
            which is already verified in \ref{o:Relative pinning maps conjugated by cocharacters}. Then, by property of $w_{a',u}$, we have 
            \[w_{a',u}U_{\alpha_{b',\frac{-l<b',a'>}{2}+m}}w_{a',u}^{-1}=U_{\alpha_{s_{a'}(b'),\frac{-l<b',a'>}{2}+m}}\;.\]
            Finally, by \ref{o:Relative pinning maps conjugated by cocharacters} again  
            \[a'^{\vee}(t^{-l/2})U_{\alpha_{s_{a'}(b'),\frac{-l<b',a'>}{2}+m}}a'^{\vee}(t^{-l/2})^{-1}=U_{\alpha_{s_{a'}(b'),\frac{l<s_{a'}(b'),a'>}{2}-\frac{l<b',a'>}{2}+m}}\]
            \[=U_{\alpha_{s_{a'}(b'),m-l<b',a'>}}=U_{s_{a',l}(\alpha_{b',m})}\;.\] 
        \end{enumerate}
    \end{Proof}
\end{lem}\bigskip 

\subsection{Main theorem}\label{subsec:Main theorem}

\begin{thm}[Main theorem]\label{t:ref{t:Construction of RGD system for the non-split case} without further requirement}
    Using the construction in \ref{n:Construction of Affine root group in non-split case}, we have that 
    \[\Bigl(\mathcal{G}(R)^{+}C_{\mathcal{G}}(S)(k),(U_{\alpha_{a',l}})_{\alpha_{a',l}\in \Phi},C_{\mathcal{G}}(S)(k)\Bigr)\] 
    is an RGD system of affine type. 
    \begin{Proof}
        We go through the axioms: 
        \begin{enumerate}
            \item[(RGD0)] We seek to show that $U_{\alpha}\neq \{1\}$ for all $\alpha\in \Phi$: It is shown in (DR1) of \cite[Satz: 27]{cite:Reducktive_Gruppen_AP} that $U_{a'}(k)=U_{\alpha_{a',0}}\neq \{1\}$. Then we extend this fact with equation \eqref{eq:coroot-move-afrg} 
            to see that $U_{\alpha_{a',n}}\neq \{1\}$ in general. 
            \item[(RGD1)] Consider, by the generalized Chevalley commutator formula for relative pinning maps, for $a',b'\in \;_{k}\Psi$ such that $ma'\neq -kb'$ for all $m,k\geq 1$ (notice, this is precisely equivalent to requiring $\alpha_{a',l}$ and $\alpha_{b,m}$ is prenilpotent pair, see \ref{n:A Summary of results about affine roots}), we have, 
            for all $l,m\in \mathbb{Z}$,  
            \[\Bigl[x_{a'}(ct^{-l}\otimes e^{a'}_{h}),x_{b'}(dt^{-m}\otimes e^{b'}_{k})\Bigr]=\prod_{i,j>0, ia'+jb'\in \;_{k}\Psi} x_{ia'+jb'}\Bigl(N_{a',b';i,j}(ct^{-l}\otimes e^{a'}_{h},dt^{-m}\otimes e^{b'}_{k})\Bigr)\] 
            \[=\prod_{i,j>0, ia'+jb'\in \;_{k}\Psi} x_{ia'+jb'}\Bigl((c^{i}d^{j})(t^{-il-jm})\otimes N_{a',b';i,j}(e^{a'}_{h},e^{b'}_{k})\Bigr)\;.\] 
            We may write $\displaystyle N_{a',b';i,j}(e^{a'}_{h},e^{b'}_{k})=\sum_{s=1}^{_{k}m_{ia'+jb'}} f_{s} e^{ia'+jb'}_{s}$ for $f_{s}\in k$ (as $V_{ia'+jb'}$ is a $k$-module), and continue the calculation above,  
            {\small\[\Bigl[x_{a'}(ct^{-l}\otimes e^{a'}_{h}),x_{b'}(dt^{-m}\otimes e^{b'}_{k})\Bigr]=\prod_{i,j>0, ia'+jb'\in \;_{k}\Psi} x_{ia'+jb'}\Bigl((c^{i}d^{j})(t^{-il-jm})\otimes (\sum_{s=1}^{_{k}m_{ia'+jb'}} f_{s} e^{ia'+jb'}_{s})\Bigr)\]} 
            \[=\prod_{i,j>0, ia'+jb'\in \;_{k}\Psi} x_{ia'+jb'}\Bigl(\sum_{s=1}^{_{k}m_{ia'+jb'}} f_{s}(c^{i}d^{j})(t^{-il-jm})\otimes  e^{ia'+jb'}_{s}\Bigr)\] 
            \[\in \Bigl\langle U_{\beta} \Big| \beta\in \{\alpha_{pa'+qb',pl+qm}\in\Phi|p,q\in\mathbb{Z}_{>0}\} \Bigr\rangle \;.\] 
            In our context, more specifically, above precisely shows that, for $(a')=\{a',2a'\}$, 
            \[\Bigl[x_{a'}(ct^{-l}\otimes e^{a'}_{i}),x_{b'}(dt^{-m}\otimes e^{b'}_{j})\Bigr]\in \Bigl\langle U_{\beta} \Big|\beta\in \{\alpha_{pa'+qb',pl+qm}\in\Phi|p,q\in\mathbb{Z}_{>0}\} \Bigr\rangle\;, \] 
            \[\Bigl[x_{a'}(ct^{-l}\otimes e^{a'}_{i}),x_{2b'}(dt^{-2m}\otimes e^{2b'}_{j})\Bigr]\in \Bigl\langle U_{\beta} \Big|\beta\in \{\alpha_{pa'+q2b',pl+q2m}\in\Phi|p,q\in\mathbb{Z}_{>0}\} \Bigr\rangle\]
            \[\subset \Bigl\langle U_{\beta} \Big|\beta\in \{\alpha_{pa'+qb',pl+qm}\in\Phi|p,q\in\mathbb{Z}_{>0}\} \Bigr\rangle\;, { and }\] 
            \[\Bigl[x_{2a'}(ct^{-2l}\otimes e^{2a'}_{i}),x_{2b'}(dt^{-2m}\otimes e^{2b'}_{j})\Bigr]\in \Bigl\langle U_{\beta} \Big|\beta\in \{\alpha_{p2a'+q2b',p2l+q2m}\in\Phi|p,q\in\mathbb{Z}_{>0}\} \Bigr\rangle\]
            \[\subset \Bigl\langle U_{\beta} \Big|\beta\in \{\alpha_{pa'+qb',pl+qm}\in\Phi|p,q\in\mathbb{Z}_{>0}\} \Bigr\rangle\;.\] 
            Recall the following elementary fact regarding commutator bracket: 
            \[[x,zy]=[x,y][x,z][[x,z],y]\text{, }[xz,y]=[x,y][[x,y],z][z,y]\text{, and }[x,y]^{-1}=[y,x]\;.\] 
            Under the construction of affine root group by generators, we, then, extend above result to general elements of the affine roots groups with the elementary facts, and obtain 
            \[\Bigl[U_{\alpha_{a',l}},U_{\alpha_{b',m}}\Bigr]\subset \Bigl\langle U_{\beta} \Big|\beta\in \{\alpha_{pa'+qb',pl+qm}\in\Phi|p,q\in\mathbb{Z}_{>0}\} \Bigr\rangle\subset U_{(\alpha_{a',l},\alpha_{b',m})}\;.\] 
            \item[(RGD2)] 
            We construct the map $m$ on each $U_{\alpha_{a',l}}^{*}$ by $m:x\in U_{\alpha_{a',l}}^{*}\mapsto w_{a',x}(t^{-l})$ (the element $w_{a',x}(t^{-l})$ is constructed in \ref{n:Preparation for (RGD2) in non-split case}). 
            We check: For all $s_{i}:=s_{\alpha_{i}}$ where $0\leq i\leq n$, and $x\in U_{\alpha_{i}}\setminus \{1\}$ (recall $\alpha_{i}$ being the simple affine roots, see \ref{n:Affine roots and affine Weyl group}), 
            we have following 
            \begin{enumerate}
                \item $m(x)\in U_{s_{i}(\alpha_{i})}xU_{s_{i}(\alpha_{i})}$: 
                This is shown in \ref{n:Preparation for (RGD2) in non-split case}. 
                \item $m(x)U_{\alpha}m(x)^{-1}=U_{s_{i}(\alpha)}$ for all $\alpha\in \Phi$: 
                This is shown in \ref{n:Preparation for (RGD2) in non-split case}. 
                \item $m(x_{1})m(x_{2})^{-1}\in C_{\mathcal{G}}(S)(k)$ for all $x_{1},x_{2}\in U_{\alpha_{i}}\setminus \{1\}$: 
                This holds by construction, take elements $x_{1}=a'^{\vee}(t^{-l/2})u_{1}a'^{\vee}(t^{-l/2})^{-1}, x_{2}=a'^{\vee}(t^{-l/2})u_{2}a'^{\vee}(t^{-l/2})^{-1}\in U_{\alpha_{a',l}}$ (see \ref{o:Isomorphism between affine root groups by coroots}), we have 
                \[m(x_{1})m(x_{2})^{-1}=w_{a',x_{1}}(t^{-l})w_{a',x_{2}}(t^{-l})\]
                \[=a'^{\vee}(t^{-l/2})w_{a',u_{1}}a'^{\vee}(t^{-l/2})^{-1}a'^{\vee}(t^{-l/2})w_{a',u_{2}}^{-1}a'^{\vee}(t^{-l/2})^{-1}\]
                \[=a'^{\vee}(t^{-l/2})w_{a',u_{1}}w_{a',u_{2}}^{-1}a'^{\vee}(t^{-l/2})^{-1}
                \in a'^{\vee}(t^{-l/2})C_{\mathcal{G}}(S)(k)a'^{\vee}(t^{-l/2})^{-1}=C_{\mathcal{G}}(S)(k)\;.\] 
            \end{enumerate}
            \item[(RGD3)] In the following of this item, we write $\alpha:=\alpha_{a',l}$. The following method is very similar to which that was utilized in \cite[lem: 5]{cite:gp_act_tb}, although the reference's aim was to prove (RGD3)* instead of (RGD3). Recall that 
            a compatible ordering always exists, in which we have 
            \[_{k}\Psi_{+}\subset j(\Psi_{+})\subset \;_{k}\Psi_{+}\cup\{0\}\;.\]  
            We work under the assumption of this ordering. 
            Recall that $\hat{k}$ is a finite separable extension of k over which $T$ splits, and for each absolute root $a\in \Psi$, we may consider $x_{a}$ to be defined over $\hat{k}$: 
            \begin{enumerate}
                \item \textbf{(When $a'\in\;_{k}\Psi_{+}$ and $l\geq 0$)}: For $a\in \eta(a')$, elements $x_{a}(ct^{-l})$'s (and $x_{b}(ct^{-2l})$'s if $b\in \eta(2a')$) of $U_a$ are uni-upper-triangular matrices with all entries above diagonal in $\hat{k}[t^{-1}]$ (see \ref{n:Existence of needed faithful rational representation for the split case}). An arbitrary element of $U_{\alpha}$, as a product of elements $x_{a}(ct^{-l})$'s and $x_{b}(ct^{-2l})$'s, is also a uni-upper-triangular matrix with all entries above diagonal in $\hat{k}[t^{-1}]$ 
                (notice $c,d\in \hat{k}$, they do not necessarily stay in $k$).  
                We have $U_{\alpha}\subset \mathcal{G}(\hat{k}[t^{-1}])$. 
                
                Diving into a subcase of the current case (a): 
                If we were to require $l\geq 1$ (instead of only $l\geq 0$) and $a'\in\;_{k}\Psi_{+}$,  
                we would have that all entries strictly above diagonal of $x_{a}(ct^{-l})$ is contained in $t^{-1}\hat{k}[t^{-1}]$ if $c\neq 0$; similarly all entries strictly above diagonal of $x_{b}(ct^{-2l})$ is contained in $t^{-2}\hat{k}[t^{-1}]$. Hence, we have, if $a'\in \Psi_{+}$ and $l\geq 1$,  $x_{a'}(ct^{-l})\in \mathcal{G}(\hat{k}[t^{-1}])$ is uni-upper-triangular with entries strictly above diagonal in $t^{-1}\hat{k}[t^{-1}]$, and $x_{2a'}(ct^{-2l})$ is uni-upper-triangular with entries strictly above diagonal in $t^{-2}\hat{k}[t^{-1}]$. 
                That is an arbitrary element of $U_{\alpha}$ is uni-upper-triangular with entries strictly above diagonal in $t^{-1}\hat{k}[t^{-1}]$. 
                
                \item \textbf{(When $a'\in\;_{k}\Psi_{+}$ and $l\leq -1$)}: An arbitrary element of $U_{\alpha}$ is a uni-upper-triangular matrix and has all entries strictly above diagonal in $t\hat{k}[t]$. Also, $U_{\alpha}\subset \mathcal{G}(\hat{k}[t])$.  
                
                \item \textbf{(When $a'\in\;_{k}\Psi_{-}$ and $l\geq 1$)}: An arbitrary element of $U_{\alpha}$ is a uni-lower-triangular matrix and has all entries strictly below diagonal in $t^{-1}\hat{k}[t^{-1}]$. Also, $U_{\alpha}\subset \mathcal{G}(\hat{k}[t^{-1}])$. 
                
                \item \textbf{(When $a'\in\;_{k}\Psi_{-}$ and $l\leq 0$)}: An arbitrary element of $U_{\alpha}$ is a uni-lower-triangular matrix and has all entries strictly below diagonal in $\hat{k}[t]$. Also, $U_{\alpha}\subset \mathcal{G}(\hat{k}[t])$. 
            \end{enumerate}
            In above, (b), (c), and (d) are obtained with similar procedure as (a). 
            Consider $U_{\epsilon}\subset \mathcal{G}(\hat{k}[t^{-\epsilon}])$ where $\epsilon=\pm$ and $t^{\pm}:=t^{\pm1}$, we have that 
            \[U_{+}\cap U_{-}\subset \mathcal{G}(\hat{k}[t^{-1}])\cap \mathcal{G}(\hat{k}[t])=\mathcal{G}(\hat{k})\;,\] 
            where the equal sign is a set-theoretic result. Consider group homomorphisms $p_{\epsilon}:\mathcal{G}(\hat{k}[t^{-\epsilon}])\to \mathcal{G}(\hat{k})$ induced by $\hat{k}[t^{-\epsilon}]\to \hat{k}$ (defined by $t^{-\epsilon}\mapsto 0$) entry-wise. 
            Considering the preimage of a subgroup under a group homomorphism is a subgroup, and (a) to (d) above, we see the following: (We denote $\mathfrak{U}^{\pm}(\hat{k})$ the upper (for +) and lower (for -) uni-triangular matrices with entries in $\hat{k}$)
            \begin{enumerate}
                \item[] In the case (a), $p_{+}(U_{\alpha_{a',0}})\subset \mathfrak{U}^{+}(\hat{k})$, 
                and, hence, $U_{\alpha}$ is contained in $p_{+}^{-1}(\mathfrak{U}^{+}(\hat{k}))$. 
                Similar reasoning give us, in case (b), the image of $U_{\alpha}$ under $p_{-}$ is $Id$, and, hence, $U_{\alpha}$ is contained in $p_{-}^{-1}(\mathfrak{U}^{-}(\hat{k}))$; 
                in case (c), the image of $U_{\alpha}$ under $p_{+}$ is $Id$, and, hence $U_{\alpha}$ is contained in $p_{+}^{-1}(\mathfrak{U}^{+}(\hat{k}))$; 
                in case (d), $U_{\alpha}$ is contained in $p_{-}^{-1}(\mathfrak{U}^{-}(\hat{k}))$. 
                To sum up, we have $U_{\epsilon}\subset p_{\epsilon}^{-1}(\mathfrak{U}^{\epsilon}(\hat{k}))$. 
                This implies that  
                \[U_{\epsilon}\cap \mathcal{G}(\hat{k})\subset p_{\epsilon}^{-1}(\mathfrak{U}^{\epsilon}(\hat{k}))\cap \mathcal{G}(\hat{k})\subset \mathfrak{U}^{\epsilon}(\hat{k})\;.\] 
                Notice the middle part is just describing elements of $p_{\epsilon}^{-1}(\mathfrak{U}^{\epsilon}(\hat{k}))$ with entries in $\hat{k}$, 
                and they have to be in $\mathfrak{U}^{\pm}(\hat{k})$ because of the construction of $p_{\epsilon}$ and the definition of $\mathfrak{U}^{\epsilon}(\hat{k})$. 
            \end{enumerate}

            We, therefore, have 
            \[U_{+}\cap U_{-}= (U_{+}\cap U_{-})\cap \mathcal{G}(\hat{k}) = 
            \Bigl(U_{+}\cap \mathcal{G}(\hat{k})\Bigr)\cap \Bigl(U_{-}\cap \mathcal{G}(\hat{k})\Bigr) \subset \mathfrak{U}^{+}(\hat{k})\cap \mathfrak{U}^{-}(\hat{k})=\{1\}\;.\]
            So we must have $U_{+}\cap U_{-\alpha_{i}}\subset U_{+}\cap U_{-}=\{1\}$. 
            The containment is because we have taken all simple affine roots to be positive, and, hence, $U_{-\alpha_{i}}\subset U_{-}$. 

            Recall (RGD0), which implies that $U_{-\alpha_{i}}\neq \{1\}$, hence, we have, for all simple affine roots $\alpha_{i}$, $U_{-\alpha_{i}}\not\leq U_{+}$. 

            As a side note, in above, we kept using $\hat{k}$ because we started our analysis with absolute pinning isomorphisms which we can take to be defined over $\hat{k}$, but not necessarily over $k$. 
            We are not claiming that the affine root groups we have defined do not land inside the $k[t,t^{-1}]$-points. 
            \item[(RGD4)] See \ref{n:Important proteries of generalize unipotent elements and affine root group} item (b) and (d), we have that 
            $\Bigl\{x_{a'}(r\otimes e^{a'}_{i}) \Big| r\in R,e^{a'}_{i}\in M^k_{a'},a'\in \;_{k}\Psi\Bigr\}$ is a generating subset of 
            $\langle U_{\alpha}|\alpha\in\Phi \rangle =\mathcal{G}(R)^{+}$. Considering $C_{\mathcal{G}}(S)(k)$ normalizes each of the affine root groups, hence, the whole elementary subgroup, the statement is proven.  
            \item[(RGD5)] See \ref{n:Important proteries of generalize unipotent elements and affine root group} item (b), where this is proven. 
        \end{enumerate}
    \end{Proof}
\end{thm}

As a further comment, to construct the RGD system for the elementary subgroup $\mathcal{G}(k[t,t^{-1}])^{+}$, one replaces $C_{\mathcal{G}}(S)(k)$'s with $C_{\mathcal{G}}(S)(k)\cap \mathcal{G}(k[t,t^{-1}])^{+}$'s in above \ref{t:ref{t:Construction of RGD system for the non-split case} without further requirement}. \bigskip

\begin{remark}\label{rmk:extending main theorem}

    In general, it is already shown (for instance, in \cite[Satz:27]{cite:Reducktive_Gruppen_AP}) that we have {\small$\mathcal{G}(k)=\mathcal{G}(k)^{+}C_{\mathcal{G}}(S)(k)$}. 
    Now if we further require that the reductive group $\mathcal{G}$ in consideration be of simply-connect type (in the same sense as in \cite{cite:hom_inv_of_non_stable_K1}), and that any normal semisimple subgroup of $\mathcal{G}$ has $k$-rank at least 2; 
    then, because \cite[cor:6.2]{cite:hom_inv_of_non_stable_K1} provides a group isomorphism {\small$\mathcal{G}(k)/\mathcal{G}(k)^+ \xrightarrow{\cong} \mathcal{G}(k[t,t^{-1}])/\mathcal{G}(k[t,t^{-1}])^+$} induced by inclusion of $\mathcal{G}(k)\hookrightarrow \mathcal{G}(k[t,t^{-1}])$, we have 
    \[\mathcal{G}(k[t,t^{-1}])=\mathcal{G}(k[t,t^{-1}])^{+}\mathcal{G}(k)=\mathcal{G}(k[t,t^{-1}])^{+}C_{\mathcal{G}}(S)(k)\;.\] 
    Hence, the RGD system we provided is, in fact, for the whole group of $\mathcal{G}(k[t,t^{-1}])$ when $\mathcal{G}$ is simply connected, and any normal semisimple subgroup of $\mathcal{G}$ has $k$-rank at least 2. 
    
\end{remark}\bigskip

%% file: sections/Examples.tex
\begin{ex}[Split case]\label{ex:Split case}
    In this example, we discuss the special case where $\mathcal{G}$ is already split over $k$; in which case the relative and absolute root systems coincide, and the projection of roots becomes the identity map, making the relative pinning maps the same as the absolute pinning isomorphisms. 
    Then \ref{t:ref{t:Construction of RGD system for the non-split case} without further requirement} shows that $(\mathcal{G}(k[t,t^{-1}])^{+}T(k), (U_{\alpha})_{\alpha\in\Phi},T(k))$, where $U_{\alpha_{a,l}}:=\{x_{a}(ct^{-l})|c\in k\}$, is an RGD system. One can see that this is essentially the case of $\mathcal{G}$ being a Chevalley group, and a proof in which case can be found at \cite[lem: 5]{cite:gp_act_tb}. 
\end{ex}\bigskip

\begin{ex}[Special Unitary group case]\label{ex:Unitary group case}
    We focus only on the case where the special unitary group has a relative root system of type BC, and we only consider fields with characteristic non 2. 
    We focus on type BC only because it is the only irreducible non-reduced root system. 
    In general, the setup is as following: 
    We take D to be a division algebra over its center $k':=Z(D)$. 
    Let D be equipped with an involution (of second type) denoted by $\tau$, and let $k$ be the fixed field of $k'$ under $\tau|_{k'}$ (hence, $D$ is also a division algebra over $k$, and can be embedded by regular representation into $M_{2n^{2}}(k)$, where $n$ in this context is the index of $D$ over $k'$). 
    Let $\epsilon:=-1$, consider the $(\epsilon,\tau)$-hermitian form denoted by $f:V\times V\to D$ on an $m$ dimensional (right) free module $V$ over $D$. 
    As introduced in \cite[sec:2.3.3]{cite:Alg_Geo_Num_Thm}, in our case, the special unitary group $\mathcal{G}(k[t,t^{-1}]):=\mathbf{SU}(V^{k(t)},f)(k[t,t^{-1}])$ can be embedded into $SL_{2mn^{2}}(\overline{k(t)})$ 
    (where $m$ is the dimension of $V$ over $D$, which is the same as the dimension of $V^{k(t)}:=V\otimes_{D}D^{k(t)}$ over $D^{k(t)}:=D\otimes_{k}k(t)$). 
    The relative root system $_{k}\Psi$ of $\mathcal{G}$ is of type $BC_{l}$ (where $l$ is the Witt index of $V$ over $D$, which is the same as the Witt index of $V^{k(t)}$ over $D^{k(t)}$), and the absolute root system $\Psi$ of $\mathcal{G}$ is of type $A_{mn-1}$. It can be seen that $\mathcal{G}$ is $\overline{k}$-isomorphic to $\mathbf{SL}_{mn}$ in this case (see case (3) of \cite[prop: 2.35]{cite:Alg_Geo_Num_Thm}). 

    We know that, for any two absolute pinning ismorphisms $x_{a}(\lambda)$ and $y_{a}(\lambda)$ for a fixed absolute root $a$, $x_{a}(\lambda)=y_{a}(c\lambda)$ for some $c\in \overline{k}^{*}$ (see \cite[18.6]{cite:borel}). Recall the well-known absolute pinning ismorphisms of $\mathbf{SL}_{mn}$ by $x_{\epsilon_{i}-\epsilon_{j}}(\lambda)=E_{ij}(\lambda)$ (roots of root system of type $A$ only have the form of $\epsilon_{i}-\epsilon_{j}$), the only thing left needed for us to determine the relative pinning maps, and, hence, the RGD system is the projection of the absolute roots to the relative roots from type $A_{mn-1}$ to type $BC_{l}$. 

    We provide an example of the needed root projection by limiting our condition to $D=k'$, and, so, $n=1$: We start our construction of the projection from the absolute roots to the relative roots on the bases of the root systems. For $A_{m-1}$, we choose 
    \[\{\epsilon_{1}-\epsilon_{2},\epsilon_{2}-\epsilon_{3},\cdots,\epsilon_{m-2}-\epsilon_{m-1},\epsilon_{m-1}-\epsilon_{m}\}\;.\] 
    For $BC_{l}$, we use 
    \[\{\epsilon'_{1}-\epsilon'_{2},\epsilon'_{2}-\epsilon'_{3},\cdots,\epsilon'_{l-2}-\epsilon'_{l-1},\epsilon'_{l-1}-\epsilon'_{l},\epsilon'_{l}\}\;.\] 
    We define the projection in a ``symmetric to the middle'' way by 
    \[\begin{cases}
        \epsilon_{i}-\epsilon_{i+1}, \epsilon_{m-i}-\epsilon_{m-i+1}\mapsto \epsilon'_{i}-\epsilon'_{i+1} & i\in [1,l-1] \\ 
        \epsilon_{l}-\epsilon_{l+1},\epsilon_{m-l}-\epsilon_{m-l+1}\mapsto \epsilon'_{l} \\ 
        \epsilon_{i}-\epsilon_{i+1}\mapsto 0 & i\in [l+1,m-l-1] 
    \end{cases}\;.\] 
    This construction generalizes to (we are only presenting the positive roots in below, adding negative sign in front all the absolute and relative roots in the following formulae provide the projections for the negative roots) 
    \[\begin{cases}
        \epsilon_{i}-\epsilon_{m-i+1}\mapsto 2\epsilon'_{i} & i\in [1,l] \\ 
        \epsilon_{i}-\epsilon_{m-j+1},\epsilon_{j}-\epsilon_{m-i+1}\mapsto \epsilon'_{i}+\epsilon'_{j} & i\neq j\in [1,l] \\ 
        \epsilon_{i}-\epsilon_{j},\epsilon_{m-j+1}-\epsilon_{m-i+1}\mapsto \epsilon'_{i}-\epsilon'_{j} & 1\leq i<j\leq l \\ 
        \epsilon_{i}-\epsilon_{h} \mapsto \epsilon'_{i} & i\in [1,l] \text{ and } h\in [l+1,m-l] \\ 
        \epsilon_{h}-\epsilon_{m-i+1} \mapsto \epsilon'_{i} & i\in [1,l] \text{ and } h\in [l+1,m-l]
    \end{cases}\;.\] 
    This results in relative pinning maps (by cases) in the following exhausted list of forms of (positve) relative roots in the relative root system of type $BC_{l}$ (the cases for negative relative roots are similar): 
    \begin{enumerate}
        \item For ``double roots'' in form of $2\epsilon'_{i}$:  
        \[\displaystyle x_{2\epsilon'_{i}}(\mu e^{2\epsilon'_{i}}_{\epsilon_{i}-\epsilon_{m-i+1}})=x_{\epsilon_{i}-\epsilon_{m-i+1}}(\mu)\;.\] 
        \item For ``addition roots'' in form of $a'=\epsilon'_{i}+\epsilon'_{j}$: 
        \[\displaystyle x_{a'}(ae^{a'}_{\epsilon_{i}-\epsilon_{m-j+1}}+be^{a'}_{\epsilon_{j}-\epsilon_{m-i+1}})=x_{\epsilon_{i}-\epsilon_{m-j+1}}(a)x_{\epsilon_{j}-\epsilon_{m-i+1}}(b)\;.\] 
        \item For ``subtraction roots'' in form of $a'=\epsilon'_{i}-\epsilon'_{j}$: 
        \[x_{a'}(ae^{a'}_{\epsilon_{i}-\epsilon_{j}}+be^{a'}_{\epsilon_{m-j+1}-\epsilon_{m-i+1}})=x_{\epsilon_{i}-\epsilon_{j}}(a)x_{\epsilon_{m-j+1}-\epsilon_{m-i+1}}(b)\;.\] 
        \item For ``single roots'' in form of $\epsilon'_{i}$: 
        \[x_{\epsilon'_{i}}\Bigl(\sum_{h=l+1}^{m-l} a_{h}e^{\epsilon'_{i}}_{\epsilon_{i}-\epsilon_{h}}+b_{h}e^{\epsilon'_{i}}_{\epsilon_{h}-\epsilon_{m-i+1}}\Bigr)=\prod^{m-l}_{h=l+1}x_{\epsilon_{i}-\epsilon_{h}}(a_{h})x_{\epsilon_{h}-\epsilon_{m-i+1}}(b_{h})\;.\]
    \end{enumerate}
    Now, one may construct the affine root groups according to \ref{n:Construction of Affine root group in non-split case} or \ref{n:Alternative form of affine root group in non-splt case}, and receive an RGD system of affine type according to \ref{t:ref{t:Construction of RGD system for the non-split case} without further requirement}. 

    As a further remark, one can also use the constructions given in \cite[10.1.2]{cite:gp_rd_sur_corp_local} to construct affine root groups as the following (we take $q$, $D_{\tau,\epsilon}$, and $Z$ from context of \cite[10.1.2]{cite:gp_rd_sur_corp_local}): 
    \[U_{\alpha_{\epsilon'_{i}\pm\epsilon'_{j},r}}:=\Bigl\{u_{i,\pm j}(\lambda t^{-r})\Big|\lambda\in D\Bigr\}\;,\] 
    \[U_{\alpha_{\epsilon'_{i},r}}:=\Bigl\{u_{i}(zt^{-r},ct^{-2r})\Big| z\in V_{0}, c\in q(z)\Bigr\}=\Bigl\{u_{i}(zt^{-r},ct^{-2r})\Big| (z,c)\in Z\Bigr\}\;,\text{ and }\] 
    \[U_{\alpha_{2\epsilon'_{i},r}}:=\Bigl\{u_{i}(0,ct^{-r})\Big| c\in D_{\tau,\epsilon}\Bigr\}\;.\] 
    One can then see, with extensive calculations according to relations of elements constructed in \cite[10.1.2]{cite:gp_rd_sur_corp_local} (some of which are listed in \cite[10.1.11]{cite:gp_rd_sur_corp_local} and \cite[pp. 108-110]{cite:twbd_and_app_to_S_A_G}), this alternative construction of affine root groups also provides an RGD system of affine type in the special unitary group case. The author conjectures that the above construction according to \cite[10.1.2]{cite:gp_rd_sur_corp_local} and the construction according to \ref{t:ref{t:Construction of RGD system for the non-split case} without further requirement} are exactly the same (in the sense that the affine root groups are exactly the same both constructions) under certain choice of absolute pinning isomorphisms. 
\end{ex}\bigskip